\documentclass{amsart}

\usepackage{tiaStyle}

\usepackage{amssymb}
\usepackage{wrapfig}
\usepackage{tikz-cd}
\usetikzlibrary{arrows}

\usepackage[yyyymmdd,hhmmss]{datetime}
\usepackage[normalem]{ulem}
\RequirePackage{color, colortbl}

\usepackage{amsopn}
\usepackage{booktabs}
\newcommand{\kron}{\otimes}
\def\be{\begin{equation}}
\def\ee{\end{equation}}
\def\bes{\begin{equation*}}
\def\ees{\end{equation*}}
\def\bea{\begin{equation} \begin{aligned}}
\def\eea{\end{aligned} \end{equation}}
\def\beas{\begin{equation*} \begin{aligned}}
\def\eeas{\end{aligned} \end{equation*}}
\def\bi{\begin{itemize}}
\def\ei{\end{itemize}}
\def\ben{\begin{enumerate}}
\def\een{\end{enumerate}}
\newtheorem{thm}{Theorem}[section]
\newtheorem{lem}{Lemma}[section]
\newtheorem{defi}{Definition}[section]

\newtheorem{prop}[thm]{Proposition}

\newtheorem{rem}{Remark}[section]
\newtheorem{cor}{Corollary}[section]
\def\bt{\begin{thm}}
\def\et{\end{thm}}
\def\bl{\begin{lem}}
\def\el{\end{lem}}
\def\bd{\begin{defi}}
\def\ed{\end{defi}}
\def\bc{\begin{cor}}
\def\ec{\end{cor}}
\def\bp{\begin{proof}}
\def\ep{\end{proof}}
\def\br{\begin{rem}}
\def\er{\end{rem}}
\newenvironment{Proof}[1][Proof]{\begin{trivlist} \item[\hskip \labelsep {\it #1.}]}{\end{trivlist}}
\def\bpp{\begin{Proof}}
\def\epp{\end{Proof}}

\RequirePackage{tikz}
  \RequirePackage{pgfplots}
    \pgfplotsset{compat=1.16}
    \pgfdeclarelayer{background}
    \pgfdeclarelayer{foreground}
    \pgfsetlayers{background,main,foreground}

\newcommand{\vertiii}[1]{{\left\vert\kern-0.25ex\left\vert\kern-0.25ex\left\vert #1
    \right\vert\kern-0.25ex\right\vert\kern-0.25ex\right\vert}}

\newcommand{\DF}{{\rm D}\!F}
\newcommand{\DtF}{{\rm D}^2\!F}

\newcommand\figurescale{1}

\title[Least-Squares Finite Element Method for ODEs]{Least-Squares Finite Element Method\\ for Ordinary Differential Equations}

\author[M.~Chung]{Matthias Chung}
\address[MC]{Department of Mathematics, Virginia Tech, Blacksburg, VA, and\\
	Academy of Data Science, Virginia Tech, Blacksburg, VA}
\email{mcchung@vt.edu}

\author[J.~Krueger]{Justin Krueger}
\address[JK]{Department of Mathematics, Virginia Tech, Blacksburg, VA}
\email{kruegej2@vt.edu}

\author[H.~Liu]{Honghu Liu}
\address[HL]{Department of Mathematics, Virginia Tech, Blacksburg, VA}
\email{hhliu@vt.edu}

\begin{document}

\begin{abstract}
We consider the least-squares finite element method (\texttt{lsfem}) for systems of nonlinear ordinary differential equations, and establish an optimal error estimate for this method when piecewise linear elements are used. The main assumptions are that the vector field is sufficiently smooth and that the local Lipschitz constant as well as the operator norm of the Jacobian matrix associated with the nonlinearity are sufficiently small, when restricted to a suitable neighborhood of the true solution for the considered initial value problem. This theoretic optimality is further illustrated numerically, along with evidence of possible extension to higher-order basis elements. Examples are also presented to show the advantages of \texttt{lsfem} compared with finite difference methods in various scenarios. Suitable modifications for adaptive time-stepping are discussed as well.
\end{abstract}


\maketitle

\section{Introduction} \label{sec:intro}
In scientific fields ranging from systems biology and systems engineering to social sciences, physical systems and finance, differential equations are omnipresent and constitute an essential tool to simulate, analyze, predict, and to ultimately make informed decisions. Due to the wide range of applications, the search for efficient, flexible, and reliable numerical schemes is still a timely topic despite its long history. Numerical solutions of ordinary differential equations (ODEs), in particular for initial-value problems (IVPs), are predominantly obtained by a rich variety of finite difference single/multistep schemes, which lead to both implicit and explicit solvers that are now standard in many programming languages \cite{ascher1994numerical,hairer2006geometric,Hairer1993,Hairer1996,lambert1991numerical}. In contrast, finite element methods for ODEs are much less investigated, despite the works on continuous and discontinuous Galerkin methods (see \cite[Section 2.2]{Cockburn00} for a brief review) and collocation methods \cite{Fairweather89}. The same can be said for delay differential equations (DDEs) as well \cite{bellen2013numerical}.

\begin{figure}[htb!]
  \begin{center}
    \includegraphics[width=0.49\textwidth]{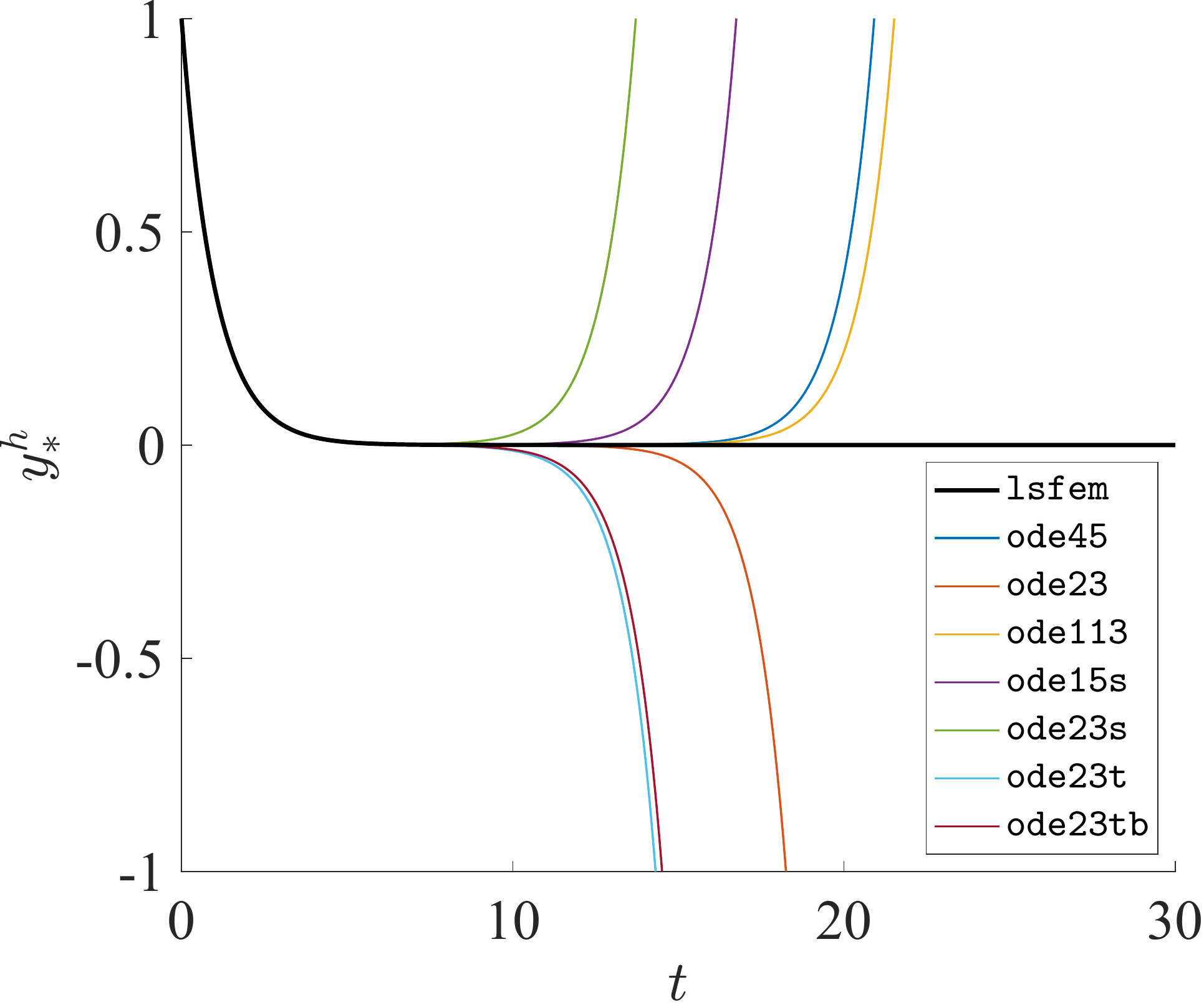}
    \includegraphics[width=0.49\textwidth]{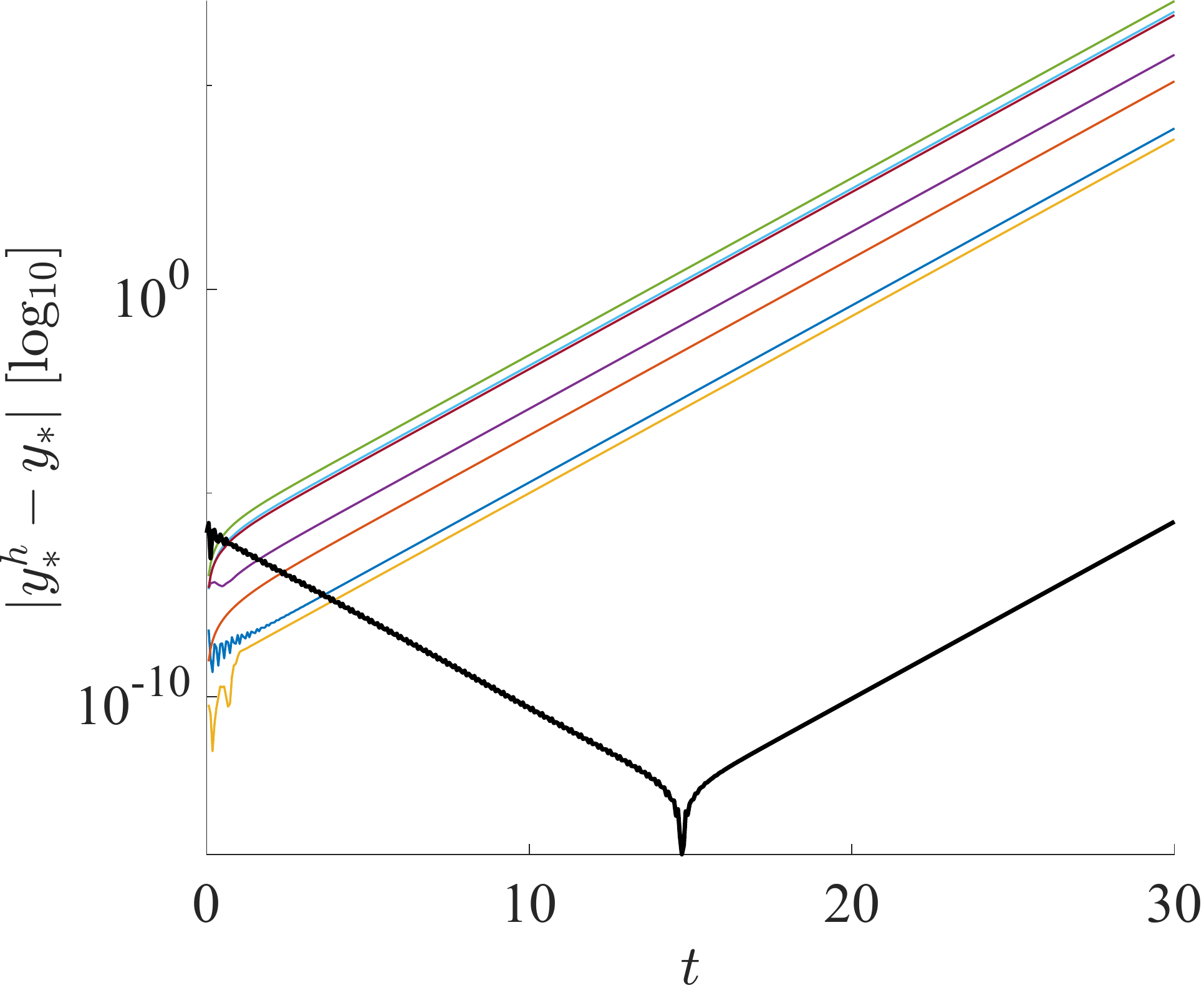}
    \caption{{\bf Failure of standard finite difference methods}. Left panel: Numerical solutions of $y' = y - 2\e^{-t}$ with $y(0) = 1$ using {\tt lsfem} solver (black bold line) and various standard ODE solvers from Matlab's ODE solver library. The solution of the IVP is denoted by $y_*$ while approximations are donoted by $y_*^h$. Right panel: The point-wise error between the exact solution $y(t) = \e^{-t}$ and the numerical solutions obtained from the numerical solvers used in the left panel. For all the Matlab's built-in ODE solvers, the relative and absolute tolerances are set to be $\texttt{RelTol} = 10^{-8}$ and $\texttt{AbsTol} =10^{-8}$, respectively; and the largest allowed step size is $\texttt{MaxStep}= 0.1$. For {\tt lsfem}, we used cubic splines as the basis functions defined on a uniform mesh of size $\delta t =0.1$.}
        \label{fig:Example1}
    \end{center}
\end{figure}

\medskip
\paragraph{\it Motivation} In this work, we initiate an effort to explore the least-squares finite element method ({\tt lsfem}) as a viable way to numerically solve ODEs and DDEs. Before entering into details, we briefly illustrate the strength of the {\tt lsfem} using a simple-looking ODE, which turns out to be challenging for traditional finite difference methods (FDMs) to operate. The problem of concern here is the following linear IVP
\begin{equation}
    y' = y - 2\e^{-t}, \quad y(0) = 1.
\end{equation}
Note that the exact solution is given by $y(t) = \e^{-t}$. Remarkably, all  Matlab built-in numerical ODE solvers fail on this example when solution over a relatively long time interval is computed. Whereas, the proposed {\tt lsfem} tracks well the exact solution. In Fig.~\ref{fig:Example1}, we present the numerical solutions (left panel) and the corresponding pointwise errors (right panel) on the interval $t \in[0,30]$ for all these solvers. As can be seen in the left panel, sooner or later, the solutions from the finite-difference schemes exhibit exponential growth, leading thus to exponentially growing pointwise errors. In contrast, the maximum error for {\tt lsfem} over the whole interval remains below $2\cdot10^{-6}$ (see black curve in the right panel of Fig.~\ref{fig:Example1}).
It is also worth noting that the setup of the experiment is actually in favor of the built-in solvers, since we used a uniform mesh size for {\tt lsfem} while allowing the build-in Matlab solvers to exhibit smaller or equal step sizes compared with the mesh size for {\tt lsfem}; see the caption of Fig.~\ref{fig:Example1} for further details.

The failure of the FDMs for the above example is actually not surprising. It results from discretization errors which are amplified exponentially over time since the equation has no stabilizing nonlinear terms to counterbalance the linear instability. Indeed, assume that at a given time instant $s>0$, the true solution $y(s) = \e^{-s}$ is perturbed by a small amount $\epsilon$, that is $\widehat{y}(s) = \e^{-s} + \epsilon$. Then, by direct calculation using the original equation, one sees that this deviation gets amplified to $\e^{t-s} \epsilon$ for all $t \ge s$.
Since local discretization errors are intrinsic to any FDMs, such deviations are unavoidable.

In contrast to the ``localization'' nature of FDMs, the aim of an {\tt lsfem} is to find an optimal approximate solution within a given subspace that minimizes an objective function over the whole time interval of integration (cf.~Section~\ref{Sect_theory_prelim}), hence making such methods much more robust to local discretization errors compared to FDMs. The {\tt lsfem} methods are also flexible in the sense that minor changes are needed when considering different types of dynamical systems, either governed by ODEs or DDEs and in the contexts of either IVPs or boundary value problems (BVPs), which allows for a unified numerical implementation for all the cases. In fact, the setup can also handle a broad class of differential algebraic equations (DAEs) as well, with the associated optimization problems become now constrained optimizations. Moreover, since the objective function directly controls the discretization error, it can be used as a diagnostic tool for local mesh adaptivity consideration, a feature crucial for problems involving abrupt local changes or stiffness.

Since their emergence in the early 1950s, finite element methods (FEMs) have become one
of the most versatile and powerful methodologies for the numerical solution of partial differential equations (PDEs). Whereas, for ODEs and DDEs, the usage of FEMs is much less pursued as mentioned above. Intuitively, this may be related to the facts that the salient feature of geometrical flexibility of FEMs is dormant in these cases, and that solutions for ODEs and DDEs are oftentimes smooth, rendering the weak formulation of FEMs less attractive.

However, as already illustrated in Fig.~\ref{fig:Example1} above, the {\tt lsfem} can provide accurate solutions in situations that traditional FDMs may fail drastically. This is further supported by other examples in Section~\ref{Sect_numerics} that the superior performance of the {\tt lsfem} reported in Fig.~\ref{fig:Example1} is not just an exception. These numerical results prompt us to re-evaluate the aforementioned intuition about the usage of FEMs, at least in the least-squares settings, for ODEs and DDEs.

These investigations are further driven by newly discovered connections between ordinary differential equations and residual neural networks \cite{haber2017stable}. Within the field of neural networks where stability is a major concern, recent works are starting to investigate finite element type solvers \cite{gunther2021spline}.

The existing literature on {\tt lsfem} is mainly devoted to PDEs; see e.g.,~\cite{Becker_al1981,Bochev2009,Davis1984,Jiang1998} and references therein. On the theoretic side, for linear PDE problems, very satisfactory theoretical understandings have already been gained that includes convergence results and even optimal error estimates \cite{Bochev2009,Jiang1998}. Nevertheless, error estimates in the case of nonlinear PDE problems remain largely open. In contrast, {\tt lsfem} for ODEs and DDEs has not received much attention yet, neither theoretically nor computationally. In this article, we take a first step in establishing {\tt lsfem} error analysis for general nonlinear ODEs, deferring the treatments for DAEs and DDEs to future works.

\medskip
\paragraph{\it Main contributions} In that respect, we consider IVP of nonlinear ODEs for which we establish under suitable conditions optimal error estimates for {\tt lsfem} with piecewise linear elements; see Theorem~\ref{Thm_errest} below. The optimality of the estimate is in the sense that the error bound established in Theorem~\ref{Thm_errest} is of the same order, in terms of the mesh size $h$, as the finite element interpolation error recalled in Lemma~\ref{Lem_projection_linearFE}. Our main idea centers around an estimate given by Proposition~\ref{Prop_errest_part1} associated with an auxiliary system \eqref{Eq_nonlinear_modified}. Given an {\tt lsfem} solution $y_*^h$ in a finite dimensional subspace $X^h$ of the Sobolev (state) space $X = H^1(0,T; \mathbb{R}^d)$, this latter auxiliary system is obtained by replacing the original nonlinear vector field $F(y) + f(t)$ in \eqref{Eq_nonlinear} by $F(y_*^h(t)) + f(t)$. Since the latter vector field consists simply of a given time-dependent function for a given $y_*^h$, optimal error estimates between the true solution $w_*$ and the {\tt lsfem} solution $w^h_*$ of the auxiliary system \eqref{Eq_nonlinear_modified} is well known using the classical Aubin-Nitsche trick \cite{aubin1967behavior,ciarlet2002finite,nitsche1968kriterium}; see Lemma~\ref{Lem_errest_F=0} and the estimate given by \eqref{Eq_nonlinear_modified_errest}, in which the dependence on $y_*^h$ are marked via $w_*[y_*^h]$ and $w^h_*[y_*^h]$. However, to establish a suitable control of the difference $\|y_*^h - w^h_*[y_*^h]\|$ between the {\tt lsfem} solutions $y_*^h$ (for the original nonlinear IVP \eqref{Eq_nonlinear}) and $w^h_*[y_*^h]$ (for the auxiliary system \eqref{Eq_nonlinear_modified}) as given in Proposition~\ref{Prop_errest_part1} requires a major effort.

Such an estimate for $\|y_*^h - w^h_*[y_*^h]\|$ is established through a series of lemmas that exploit geometric properties revealed by the first-order optimality condition associated with each minimizer $y_*^h$ in the subspace $X^h$ for the objective function $J$ given by \eqref{Def_J}. Indeed, from $\frac{\d}{\d \tau} J(y_*^h+\tau v_h; F, f, g) \big\vert_{\tau = 0} = 0$ for all $v_h$ in $X^h$, after some algebraic operations, we can actually link this necessary condition with $w^*[y_*^h]$ through the following orthogonality property (cf.~Lemma~\ref{Lem_geometric_property}):
\be \label{Eq_EL_identity_intro}
\langle y_*^h - w_*[y_*^h], \; v_h - \Gamma(\cdot; v_h) \rangle_X = 0, \quad \forall \; v_h \in X^h,
\ee
where $\Gamma$ is an integral involving the Jacobian matrix of $F$ given by \eqref{Eq_Gamma}. It is this simple, albeit not so obvious, geometric identity that opens the room for estimation, once $\Gamma$ is further split as the sum of its projection $\Pi_h \Gamma$  onto $X^h$ and its orthogonal complement $\Pi^\perp_h \Gamma$; see Lemmas~\ref{Lem_EL_identity_RHS} and \ref{Lem_EL_identity_LHS}.

Although the error analysis presented in this article focuses on piecewise linear elements, numerical evidence provided in Section~\ref{Sect_numerics} indicates that when a piecewise spline basis of degree $k$ is used to form $X^h$ and $F$ is $C^{k+1}$-smooth, then the error bound scales like $h^{k+1}$. Rigorous justification of such an error estimate will be addressed in a future work.

\medskip
\paragraph{\it Organization} This article is organized as follows. We first recall in Section~\ref{sec:theory} the basic setup of {\tt lsfem} in the context of IVP for nonlinear ODE systems. Besides its functional framework recalled in subsection~\ref{Sect_theory_prelim}, for later usage we also present in subsection~\ref{Sect_theory_conv} a result concerning the convergence of {\tt lsfem} solutions to the true solution; see Theorem~\ref{Thm_convergence}. While the treatment makes a direct usage of a general convergence result on the approximation of abstract nonlinear equations (cf.~\cite[Theorem 3.3, p.307]{girault2012finite} and \cite[Theorem 8.1]{Bochev2009}) some detailed calculation is required to recast the problem into the functional form dealt with in \cite[Theorem 3.3, p.307]{girault2012finite} and also to check the required assumptions therein. We provide thus a proof of this convergence theorem in~\ref{Sect_conv_thm_proof} for the sake of clarity. The associated optimal error analysis reviewed above is then dealt with in Section~\ref{Sect_theory_err_estimates}. The algorithmic aspects are then presented in Section~\ref{Sect_numerics} (cf.~Algorithm~\ref{alg:lsfem}) along with numerical results for various concrete examples that confirm the error bounds obtained in Section~\ref{Sect_theory_err_estimates} and also provide numerical evidence for possible extension to higher-order basis elements. We also discuss within this section suitable modifications for adaptive time stepping. Finally, Section~\ref{Sect_conclusion} provides a brief conclusion and potential future directions.

\section{Preliminaries} \label{sec:theory}

As a preparation for later sections concerning the error estimates (Section~\ref{Sect_theory_err_estimates}) as well as the numerical treatments (Section~\ref{Sect_numerics}), we briefly summarize the basic setup for {\tt lsfem} of first-order ODEs and then recall a classical convergence result for the {\tt lsfem} solutions.
For ease of reference, a table of the main symbols used in this work is provided in Table~\ref{tab_symbol_list}.

\begin{table}[h]
\caption{List of main symbols}
\label{tab_symbol_list}
\centering
{\footnotesize
\begin{tabular}{ll}
\toprule\noalign{\smallskip}
$X$ & The Sobolev space $H^1(0,T; \mathbb{R}^d)$ equipped with  the inner product \eqref{Eq_X_inner} and \\
      &   the corresponding induced norm \eqref{Eq_X_norm} \vspace{0.25em}  \\
$X^h$ & A finite element subspace of $X$  \vspace{0.25em}  \\
$I_h$ & The interpolation operator from $X$ to $X^h$   \vspace{0.25em} \\
$\Pi_h$ & The orthogonal projection from $X$ to $X^h$  \vspace{0.25em} \\
$\mathrm{Id}_X$ & The identity map on $X$  \vspace{0.25em} \\
$\Pi^\perp_h$ & Orthogonal complement of $\Pi_h$: $\Pi^\perp_h = \mathrm{Id}_X - \Pi_h$   \vspace{0.25em} \\
$y_*$ & Solution to the variational formulation \eqref{Eq_opt} of the IVP \eqref{Eq_nonlinear}    \vspace{0.25em} \\
$y^h_*$ & {\tt lsfem} approximation of $y_*$ in the subspace $X^h$; i.e., solution of \eqref{Eq_opt_FEM}   \vspace{0.25em} \\
$w_*[y^h_*]$ & Solution of the auxiliary system \eqref{Eq_nonlinear_modified}   \vspace{0.25em} \\
$w^h_*[y^h_*]$ & {\tt lsfem} approximation of $w_*[y^h_*]$ in the subspace $X^h$ \vspace{0.25em} \\
$u$  & A generic element in $X$ or the solution of \eqref{Eq_1st_order_cond_operator_form} depending on  the context \vspace{0.25em}  \\
$v$ & A generic element in $X$  \vspace{0.25em} \\
$\mathfrak{B}$ & The subset in $\mathbb{R}^d$ defined by \eqref{Def_set_B}, which contains both $y_*(t)$ and \\
   			& the {\tt lsfem} solution $y_*^h(t)$ for all $t$ in $[0,T]$ and all sufficiently small $h$  \vspace{0.25em} \\
$\mathfrak{C}$ & Embedding constant for the continuous embedding from $X$ to $C([0,T]; \mathbb{R}^d)$  \vspace{0.25em} \\
$\widetilde{\mathfrak{C}}$ & Embedding constant for the continuous embedding from $X $ to $L^2(0,T; \mathbb{R}^d)$  \vspace{0.25em} \\
$\langle \cdot, \cdot \rangle$ &  The standard dot product on $\mathbb{R}^d$    \vspace{0.25em} \\
$\|\cdot\|$ &  The Euclidean norm on $\mathbb{R}^d$   \vspace{0.25em} \\
$\|\cdot\|_{\mathrm{op}}$ &  The operator norm for a $d\times d$ matrix, i.e., $\|M\|_{\mathrm{op}} = \sup_{z \in \mathbb{R}^d, \|z\| = 1} \|M z\|$ \vspace{0.25em} \\
$L(Y,Z)$ & The set of bounded linear maps from a Hilbert space $Y$ to a Hilbert space $Z$  \\
\noalign{\smallskip} \bottomrule
\end{tabular}}
\end{table}

\subsection{Formulation of {\tt lsfem}} \label{Sect_theory_prelim}
We provide in this subsection a brief account of the {\tt lsfem} for first-order (nonlinear) ODE systems; and refer to \cite[Chap.~3]{Jiang1998} for more details. Given a fixed $T > 0$, consider the following initial-value problem (IVP) in $\mathbb{R}^d$ for some $d \in \mathbb{N}$:
\bea\label{Eq_nonlinear}
& y' = F(y) + f(t), && t \in (0, T],  \\
& y(0) = g,
\eea
where $F\colon \mathbb{R}^d \rightarrow \mathbb{R}^d$ is a given smooth and possibly nonlinear function, $f$ is a function in $L^2(0, T; \mathbb{R}^d)$, and $g$ is a given vector in $\mathbb{R}^d$. Precise smoothness on $F$ will be specified later on, and additional regularity on $f$ will be added when optimal error estimates are considered in Section~\ref{Sect_theory_err_estimates}.

Before proceeding, it is worth mentioning that all the results of both the current section and Section~\ref{Sect_theory_err_estimates} hold for more general systems of the form $y' = F(t,y) + f(t)$ as well. Indeed, by introducing an auxiliary scalar equation $p' = 1$ supplemented with $p(0) = 0$, and considering the new variable $z = (p,y)^\top$, we get $z' = \widetilde{F}(z) + \widetilde{f}(t)$ with $\widetilde{F}(z) =  [1, F(p,y)]^\top$ and $\widetilde{f} = [0, f]^\top$. This latter system for $z$ is an equivalent formulation of the original problem and fits into the form given by \eqref{Eq_nonlinear}.

Throughout the article, we denote the classical Sobolev space $H^1(0,T; \mathbb{R}^d)$ by $X$, which consists of $L^2(0,T; \mathbb{R}^d)$ functions whose first-order weak derivative is also in $L^2(0,T; \mathbb{R}^d)$. $X$ will be equipped with a norm that is equivalent to the usual $H^1$-norm; see \eqref{Eq_X_norm} below. Recall that a function $y \in X$ is called a strong solution of \eqref{Eq_nonlinear} if $y(0) = g$, and $y' = F(y) + f(t)$ for almost every $t \in (0 ,T)$.

The {\tt lsfem} for the IVP \eqref{Eq_nonlinear} relies on a variational reformulation of the ODE system, which seeks for $y_* \in X$ that minimizes the following objective function
\be \label{Def_J}
J(y; F, f, g) := \tfrac{1}{2}\|y' - F(y) - f\|^2_{L^2(0, T; \mathbb{R}^d)} +\tfrac{1}{2}\|y(0) - g\|^2,
\ee
where $\|\cdot\|$ denotes the Euclidean norm on $\mathbb{R}^d$. Note that if the IVP \eqref{Eq_nonlinear} admits a unique strong solution in $X$, then this solution is also the unique solution of the following unconstrained minimization problem:
\be \label{Eq_opt}
\text{Find} \quad \argmin_{y \in X} \ J(y; F, f, g).
\ee
Given any finite element subspace $X^h$ of $X$, with $h$ denoting the maximal length of the finite elements, the {\tt lsfem} for the IVP \eqref{Eq_nonlinear} consists of solving
 the following analogue of the unconstrained minimization problem \eqref{Eq_opt} restricted to $X^h$:
\be \label{Eq_opt_FEM}
\text{Find} \quad  \argmin_{y^h \in X^h} \ J(y^h; F, f, g).
\ee

Let us introduce the following inner product on $X$, which is naturally related to the objective function $J$ defined in \eqref{Def_J}, i.e.,
\be \label{Eq_X_inner}
\langle u, v \rangle_{X} = \int_0^T \langle u'(t), v'(t) \rangle \ \d t + \langle u(0), v(0) \rangle, \quad \forall \; u,v \in X,
\ee
where $\langle \cdot, \cdot \rangle$ denotes the dot product in $\mathbb{R}^d$. The norm on $X$ induced by the above inner product $\langle \cdot, \cdot \rangle_{X}$ will be denoted by $\| \cdot \|_X$, which is often referred to as the energy norm  in the literature, namely,
\be \label{Eq_X_norm}
\|u\|_X = \sqrt{\langle u, u \rangle_{X}} = \left(\int_0^T \langle u'(t), u'(t) \rangle \ \d t + \langle u(0), u(0) \rangle \right)^{1/2}, \quad u \in X.
\ee
One can check by using basic Sobolev inequalities that the $\| \cdot \|_X$-norm is equivalent to the usual Sobolev norm on $H^1(0,T; \mathbb{R}^d)$ defined by $\|u\|_{H^1} = \big(\|u\|^2_{L^2(0,T; \mathbb{R}^d)} + \|u'\|^2_{L^2(0,T; \mathbb{R}^d)}\big)^{1/2}$. Note, there exist positive constants $c_1$ and $c_2$ such that for all $u \in X$ it holds that $c_1 \|u\|_{X} \le \|u\|_{H^1} \le c_2 \|u\|_{X}$.

For later usage, let us also introduce two embedding constants. First note that since $H^1(0,T; \mathbb{R}^d)$ is continuously embedded into $C([0,T]; \mathbb{R}^d)$, see e.g.,~\cite[Theorem 8.8]{brezis_book}, then $X$ equipped with the norm defined in \eqref{Eq_X_norm} is also continuously embedded into $C([0,T]; \mathbb{R}^d)$. Throughout this article, we denote by $\mathfrak{C}$ the associated embedding constant, where $\mathfrak{C}$ is the smallest constant such that\footnote{For each $u\in X$, we always consider its continuous representative in the corresponding equivalent class. There exists a unique such representative for each $u \in X$; cf.~\cite[Theorem 8.2]{brezis_book}.}
\be \label{Eq_embedding_const}
\max_{t\in[0,T]}\|u(t)\|  \le \mathfrak{C} \|u\|_X, \quad \forall \, u \in X.
\ee
We denote also by $\widetilde{\mathfrak{C}}$ the embedding constant for the continuous embedding from $X$ to $L^2(0,T; \mathbb{R}^d)$, which is the smallest constant such that
\be \label{Eq_embedding_const2}
\|u\|_{L^2(0,T; \mathbb{R}^d)}  \le \widetilde{\mathfrak{C}} \|u\|_X, \quad \forall \, u \in X.
\ee

\subsection{Convergence of the {\tt lsfem} solutions} \label{Sect_theory_conv}
To prepare for the error analysis carried out in Section~\ref{Sect_theory_err_estimates}, we summarize in this subsection a convergence theorem for the {\tt lsfem} solutions as the dimension of the subspace $X^h$ in \eqref{Eq_opt_FEM}  increases.
The treatment makes a direct use of a general result on approximation of abstract nonlinear equations; cf.~\cite[Theorem 3.3, p.307]{girault2012finite} and \cite[Theorem 8.1]{Bochev2009}.

We work with a sequence of finite element subspaces $\{X^h \subset X\}$, with $h$ denoting the maximal length of the finite elements, such that
\be \label{Eq_Pih_conv}
\lim_{h \rightarrow 0} \|(\mathrm{Id}_X - \Pi_h ) v\|_X  = 0, \quad \forall \; v \in X,
\ee
where $\Pi_h \colon X \rightarrow X^h$ denotes the orthogonal projection onto $X^h$ under the inner product $\langle \cdot, \cdot \rangle_X$ defined in \eqref{Eq_X_inner}.

We denote by $\DF$ the Jacobian matrix of $F$, and by $\|\cdot\|_{\mathrm{op}}$ the operator norm of a bounded linear map from $\mathbb{R}^d$ onto itself.

\bt \label{Thm_convergence}
Consider the IVP \eqref{Eq_nonlinear}. Assume that $f \in L^2(0, T; \mathbb{R}^d)$, $F: \mathbb{R}^d \rightarrow \mathbb{R}^d$ is $C^3$ smooth, and \eqref{Eq_nonlinear} has a unique strong solution $y_*$ in $X$.  Assume also that $\|\DF(y_*(t))\|_{\mathrm{op}}$ is sufficiently small for all $t \in [0, T]$. Let $\mathcal{O}$ be any given open neighborhood of $y_*$ in $X$, and $\{X^h \subset X \}$ be a sequence of finite element subspaces satisfying \eqref{Eq_Pih_conv}. Then problem~\eqref{Eq_opt_FEM} has a unique solution $y_*^h$ in $\mathcal{O}$ for all sufficiently small $h$, and $y_*^h$ converges in $X$-norm to the solution $y_*$ of~\eqref{Eq_opt} as $h$ is reduced to zero,
\be \label{Eq_conv_Sect2}
\lim_{h \rightarrow 0} \|y_* - y_*^h\|_X = 0.
\ee
\et

Since some detailed calculation is required to recast the problem into the functional form dealt with in \cite[Theorem 3.3, p.307]{girault2012finite} and also to check the required assumptions therein, we provide a proof of the above theorem in \ref{Sect_conv_thm_proof} for the sake of clarity.

With the above convergence result available, we are ready to address the associated error analysis. In particular, we show for the case of piecewise linear elements that the {\tt lsfem} achieve optimal rate of convergence, which is the rate dictated by the interpolation error.


\section{Optimal {\tt lsfem} error estimates for nonlinear ODEs} \label{Sect_theory_err_estimates}

In this section, we derive an optimal error estimates for the {\tt lsfem} solutions for first-order nonlinear ODE system of the form \eqref{Eq_nonlinear}. The results are obtained for piecewise linear finite elements. Under suitable assumptions, it is shown that the error bound for {\tt lsfem} solutions is proportional to the square of the mesh size, which is of the same order as the interpolation error for piecewise linear finite elements.

Let us first introduce the following assumption about the IVP \eqref{Eq_nonlinear}:
\bi
\item[{\bf (A1)}] $f \colon [0,T] \rightarrow \mathbb{R}^d$ is absolutely continuous, $f'$ belongs to $L^2(0,T; \mathbb{R}^d)$, and $F \colon \mathbb{R}^d \rightarrow \mathbb{R}^d$ is $C^3$ smooth.  The IVP \eqref{Eq_nonlinear} has a unique solution $y_*$ in $X$.
\ei
Except the strengthened smoothness and integrability requirements on $f$, the other parts in Assumption {\bf (A1)} are the same as those required in Theorem~\ref{Thm_convergence}.

In Theorem~\ref{Thm_convergence}, a smallness assumption is also made on, $\|\DF(y_*(t))\|_{\mathrm{op}}$, the operator norm of the Jacobian matrix $\DF$ along the solution trajectory $y_*$. For the derivation of error estimates, this technical assumption needs to be further strengthened and augmented to require that both $\|\DF\|_{\mathrm{op}}$ and the local Lipschitz constant of $F$ are sufficiently small over a bounded set in $\mathbb{R}^d$ that contains the solution $y_*$ as well as the {\tt lsfem} solutions for all time $t \in [0, T]$.

We make precise these smallness assumptions on $F$ below for the sake of clarity.
Let us first note that the smallness of $\|\DF(y_*(t))\|_{\mathrm{op}}$ required in Theorem~\ref{Thm_convergence} is made precise in its proof given by \ref{Sect_conv_thm_proof}. It suffices to require that (see \eqref{eq_M_cond})
\be \label{eq_M_cond_recall}
\sup_{t \in [0,T]} \| \DF(y_*(t))\|_{\mathrm{op}} < \frac{1}{\sqrt{2T^2 + T + \widetilde{\mathfrak{C}}^2 + \sqrt{(2T^2 + T + \widetilde{\mathfrak{C}}^2)^2 + 2 T \widetilde{\mathfrak{C}}^2(1 + 2T)}}},
\ee
where $\widetilde{\mathfrak{C}}$ denotes the embedding constant for the continuous embedding from $X$ to $L^2(0,T; \mathbb{R}^d)$; cf.~\eqref{Eq_embedding_const2}.

To present the needed augmentations of \eqref{eq_M_cond_recall}, we first establish some notations which will be used throughout this section. We take the neighborhood $\mathcal{O}$ of $y_*$ in Theorem~\ref{Thm_convergence} to be an open ball in $X$ centered at $y_*$ with some radius $r>0$, which is denoted by $B(y_*,r)$. Let $\overline{h}>0$ be chosen such that for each $h \in (0, \overline{h})$, the {\tt lsfem} problem \eqref{Eq_opt_FEM} has a unique solution $y_*^h$ in $B(y_*,r)$; the existence of such an $\overline{h}$ is guaranteed by Theorem~\ref{Thm_convergence}.

With the embedding constant $\mathfrak{C}$ that ensures \eqref{Eq_embedding_const}, we define then
\be \label{Def_set_B}
\mathfrak{B} = \mathop{\bigcup}_{t \in [0, T]} \big\{ p \in \mathbb{R}^d \; :\;  \|p - y_*(t)\| < \mathfrak{C}\, r \big\}.
\ee
Since $y^h_*$ stays in $B(y_*,r) \subset X$ for all $h \in (0, \overline{h})$, it holds that $\|y^h_*(t)-y_*(t)\| \le \mathfrak{C} \|y^h_*-y_*\|_X < \mathfrak{C} \, r$, namely,
\be \label{Eq_LSFEM_in_B}
y^h_*(t) \in \mathfrak{B}, \quad \forall \, h \in (0, \overline{h}), \; t \in [0, T].
\ee

The aforementioned smallness assumptions on $F$ are as follows:
\bi
\item[{\bf (A2)}] Assume that \eqref{eq_M_cond_recall} holds. Let $r>0$ be arbitrarily given and $\overline{h}>0$ be chosen so that the {\tt lsfem} solution $y_*^h \in X^h$ stays in the ball $B(y_*,r) \subset X$ for all $h\in (0, \overline{h})$. Let $\mathfrak{B}$ be the subset in $\mathbb{R}^d$ defined by \eqref{Def_set_B} that contains both $y_*(t)$ and the {\tt lsfem} solutions $y_*^h(t)$ for all $h \in (0, \overline{h})$ and all $t\in[0,T]$; cf.~\eqref{Eq_LSFEM_in_B}. Assume that the local Lipschitz constant of $F$ over $\mathfrak{B}$ satisfies
\be \label{Eq1_in_A2}
\frac{\mathrm{Lip}(F\vert_\mathfrak{B}) T}{\sqrt{2}} < 1,
\ee
and that its Jacobian matrix satisfies
\be \label{Eq2_in_A2}
\sup_{z \in \mathfrak{B}}\|\DF(z)\|_{\mathrm{op}} < \frac{1}{\widetilde{\mathfrak{C}}},
\ee
where $\|\cdot\|_{\mathrm{op}}$ denotes the operator norm of a bounded linear map from $\mathbb{R}^d$ onto itself, and $\widetilde{\mathfrak{C}}$
is the same as given in \eqref{eq_M_cond_recall}.
\ei

Of course, all the three conditions \eqref{eq_M_cond_recall}, \eqref{Eq1_in_A2}, and \eqref{Eq2_in_A2} in {\bf (A2)} can be summarized into one assumption of the form $\sup_{z \in \mathfrak{B}}\|\DF(z)\|_{\mathrm{op}} < C$ with $C$ taken to be the right-hand-side (RHS) of \eqref{eq_M_cond_recall}. However, we prefer to keep them separate in the hope for future improvements since they are used in separate parts of the proof.

The main result of this section is summarized in the following theorem.
\bt \label{Thm_errest}
Given a sequence of subspaces $\{X^h \subset X \}$ satisfying \eqref{Eq_Pih_conv} and spanned by piecewise linear basis functions, let us consider for each $X^h$ the least-squares finite element approximation \eqref{Eq_opt_FEM} of the nonlinear IVP \eqref{Eq_nonlinear}.
Assume the assumptions {\bf (A1)} and {\bf (A2)} hold. Let $\overline{h}$ be as given in {\bf (A2)}. Then, there exists a constant $C>0$ independent of $h$ such that the {\tt lsfem} solution $y_*^h$ of \eqref{Eq_nonlinear} satisfies:
\be \label{Eq_errest_goal}
\|y_*^h - y_*\|_{L^2(0, T; \mathbb{R}^d)} \le C h^2,  \qquad \forall \; h \in (0, \overline{h}).
\ee
\et

We first recall a well known error estimate for the special case of the IVP \eqref{Eq_nonlinear}, in which $F$ is identically zero. We consider for the moment
\bea\label{Eq_trivial}
& \widetilde{y} \,' =  f(t), && t \in (0, T],  \\
& \widetilde{y}(0) = g.
\eea
In this case, its solution is obviously given by
\be \label{Eq_soln_linear}
\widetilde{y}_*(t) = g + \int_0^t f(s) \ \d s, \quad t \in [0, T].
\ee
For a given finite element subspace $X^h$, the corresponding {\tt lsfem} approximation $\widetilde{y}^h_*$ is obtained by solving
\be \label{Eq_soln_linear_LSFEM}
\argmin_{\widetilde{y}^h \in X^h}  \ \tfrac{1}{2}\|(\widetilde{y}^h)'  - f\|^2_{L^2(0, T; \mathbb{R}^d)} +\tfrac{1}{2}\|\widetilde{y}^h(0) - g\|^2.
\ee

\bl \label{Lem_errest_F=0}
Consider the problem \eqref{Eq_trivial}. Assume that $f \colon [0,T] \rightarrow \mathbb{R}^d$ is absolutely continuous and $f'$ belongs to $L^2(0,T; \mathbb{R}^d)$. Assume also that $X^h$ is spanned by piecewise linear basis functions. Then, for each such $X^h$, there exists a unique {\tt lsfem} solution that solves \eqref{Eq_soln_linear_LSFEM}, which is given by $\widetilde{y}^h_* = \Pi_h \widetilde{y}_*$, where $\widetilde{y}_*$ is the the solution of \eqref{Eq_trivial} and $\Pi_h$ denotes the orthogonal projection from $X$ onto $X^h$. Moreover, there exists a positive constant $C$ independent of $h$ such that
\be \label{Eq_error_estimate_trivial_case}
\|\widetilde{y}_* - \widetilde{y}_*^h\|_{L^2(0, T; \mathbb{R}^d)} \le C h^2 \|\widetilde{y}_*''\|_{L^2(0, T; \mathbb{R}^d)}.
\ee
\el

Although the above results are classical, we provide in \ref{Sect_Linear_case_estimation} some elements of the proof for the sake of completeness.

Note that since the $L^2$-error of piecewise linear interpolation for general function in $H^2(0, T; \mathbb{R}^d)$ is of the order $h^2$ (cf.~Lemma~\ref{Lem_projection_linearFE} below), the above result shows that the corresponding {\tt lsfem} provides the optimal convergence rate for the special case \eqref{Eq_trivial}. However, the proof of Lemma~\ref{Lem_errest_F=0} admits no straightforward extension to the general nonlinear case.

 To bridge the gap between the setting of Lemma~\ref{Lem_errest_F=0} (dealing with $F=0$) and that of Theorem~\ref{Thm_errest} (dealing with general nonlinear $F$), we introduce now an auxiliary system that will serve as a pivot in the estimates presented  below. We consider, for a given {\tt lsfem} solution $y_*^h$ of the IVP~\eqref{Eq_nonlinear}, the following auxiliary system
\bea \label{Eq_nonlinear_modified}
& w' = F(y_*^h(t)) + f(t),  && t \in (0, T], \\
& w(0) = g.
\eea

First note that the IVP \eqref{Eq_nonlinear_modified} fits into the form of \eqref{Eq_trivial} since $F(y_*^h(t)) + f(t)$ is known once $y_*^h$ is given. Lemma~\ref{Lem_errest_F=0} is thus applicable. It follows that \eqref{Eq_nonlinear_modified} always admits a unique {\tt lsfem} solution $w_*^h[y_*^h]$ under the assumption of Lemma~\ref{Lem_errest_F=0}. Moreover, denoting by $w_*[y_*^h]$ the solution of  \eqref{Eq_nonlinear_modified}, it holds that
\be \label{Eq_nonlinear_modified_FEM_soln}
w_*^h[y_*^h] = \Pi_h w_*[y_*^h],
\ee
and that
\be \label{Eq_nonlinear_modified_errest}
\|w_*[y_*^h] - w_*^h[y_*^h]\|_{L^2(0, T; \mathbb{R}^d)} \le C h^2 \|(w_*[y_*^h] )''\|_{L^2(0, T; \mathbb{R}^d)}.
\ee

We aim to derive the following estimate of $\|y_*^h - w_*^h[y_*^h]\|_{L^2(0, T; \mathbb{R}^d)}$.

\begin{prop} \label{Prop_errest_part1}
For each given $X^h$, the following estimate holds for the {\tt lsfem} solution $y_*^h$ of the nonlinear problem \eqref{Eq_nonlinear} and the {\tt lsfem} solution $w_*^h[y_*^h]$ of the problem \eqref{Eq_nonlinear_modified}:
\be \label{Eq_errest_part1_goal}
\|y_*^h -  w_*^h[y_*^h]\|_X \le \frac{C \|(w_*[y_*^h])''\|_{L^2(0, T; \mathbb{R}^d)} }{1 - \widetilde{\mathfrak{C}} \sup_{z \in \mathfrak{B}}\|\DF(z)\|_{\mathrm{op}}} h^2,
\ee
where $C>0$ denotes a universal constant independent of $h$, and $\widetilde{\mathfrak{C}}$ denotes the embedding constant between $L^2(0,T; \mathbb{R}^d)$ and $X$.
\end{prop}

We present next a few lemmas that will be used in the proof of the above Proposition.
\bl \label{Lem_projection_linearFE}
Let $X^h$ be a subspace of $X$ spanned by piecewise linear basis functions. Given any $f \in X$, denote by $I_h f$ the interpolant of $f$ in $X^h$. Then, there exists a positive constant $C$ independent of $h$ such that the following inequalities hold for all $f$ in the subspace $H^2(0,T; \mathbb{R}^d) \subset X$:
\begin{subequations} \label{Eq_projection_error_linearFE}
\begin{align}
& \|f - I_h f\|_{L^2(0,T; \mathbb{R}^d)} \le C h^2 \|f''\|_{L^2(0, T; \mathbb{R}^d)},  \label{Eq_projection_error_linearFE_a}\\
& \|f' - (I_h f)'\|_{L^2(0,T; \mathbb{R}^d)} \le C h \|f''\|_{L^2(0, T; \mathbb{R}^d)},  \label{Eq_projection_error_linearFE_b} \\
& \|f - \Pi_h f\|_{X} \le C h \|f''\|_{L^2(0, T; \mathbb{R}^d)}.  \label{Eq_projection_error_linearFE_c}
\end{align}
\end{subequations}
\el
The first two inequalities above are classical; see e.g.,~\cite[Section 2.5]{Jiang1998} for a proof. The estimate \eqref{Eq_projection_error_linearFE_c} follows from \eqref{Eq_projection_error_linearFE_b} by noting that
\be \label{Eq_projection_vs_interpolation}
\|f - \Pi_h f\|_{X} \le \|f - I_h f\|_{X} = \|f' - (I_h f)'\|_{L^2(0,T; \mathbb{R}^d)}.
\ee
The first inequality in \eqref{Eq_projection_vs_interpolation} holds because for any $f\in X$, its projection $\Pi_h f$ minimizes the residual error $\| f - v^h\|_X$ among all $v^h \in X^h$. Note also that since $t=0$ is an interpolation point of the piecewise linear finite element subspace $X^h$, it holds that $f(0) - I_h f(0) = 0$. The second equality in \eqref{Eq_projection_vs_interpolation} follows.

\bl \label{Lem_geometric_property}
Let $y_*^h$ be the {\tt lsfem} solution of the nonlinear problem \eqref{Eq_nonlinear} as specified in Theorem~\ref{Thm_errest}. Let $w_*[y_*^h]$ be the solution of the auxiliary problem \eqref{Eq_nonlinear_modified}. We define
\be \label{Eq_Gamma}
\Gamma(t; v^h) = \int_0^t \DF(y_*^h(s)) v^h(s)\ \d s, \qquad t \in [0, T], \; v^h \in X^h.
\ee
Then, the following identity holds
\be \label{Eq_EL_identity}
\langle y_*^h - w_*[y_*^h], \; v^h - \Gamma(\,\cdot\,; v^h) \rangle_X = 0, \quad \forall \; v^h \in X^h.
\ee
\el

\begin{Proof}
The equality \eqref{Eq_EL_identity} is just a reformulation of the first-order necessary condition for $y_*^h$ to be a solution of the minimization problem \eqref{Eq_opt_FEM}. Indeed, note that this latter condition is given by
\be \label{Eq_1st_order_cond_FEM_recall}
\int_0^T \langle (y^h_*)' - F(y^h_*) - f, (v^h)' - \DF(y^h_*) v^h \rangle \ \d t + \langle y^h_*(0) - g, v^h(0) \rangle = 0,  \; \forall \; v^h \in X^h,
\ee
see~\eqref{Eq_1st_order_cond_FEM} in Appendix~\ref{Sect_conv_thm_proof}.
Note also that
\bes
w_*[y_*^h] = g + \int_0^t F(y^h_*(s)) + f(s) \ \d s.
\ees
Then, \eqref{Eq_EL_identity} follows from \eqref{Eq_1st_order_cond_FEM_recall} by simply noting that $(w_*[y_*^h])'(t) = F(y^h_*(t)) + f(t)$, $\Gamma'(t; v^h) = \DF(y^h_*(t)) v^h(t)$, $w_*[y_*^h](0) =g$, and $\Gamma(0; v^h) = 0$.
\qed
\end{Proof}

The above identity \eqref{Eq_EL_identity} serves as the starting point of our estimates for the term $y_*^h - w_*^h[y_*^h]$. For this purpose, we split $\Gamma$ defined by \eqref{Eq_Gamma} as
\be \label{Eq_Gamma_split}
\Gamma(t; v^h) = \Pi_h \Gamma(t; v^h) + \Pi_h^\perp \Gamma(t; v^h),
\ee
where $\Pi_h^\perp = \mathrm{Id}_X - \Pi_h$.

To simplify the notations, we also denote
\be \label{Eq_gamma}
\gamma = y_*^h - w_*[y_*^h].
\ee

Using \eqref{Eq_Gamma_split} and \eqref{Eq_gamma} in \eqref{Eq_EL_identity}, we obtain
\be \label{Eq_EL_identity_rearranged}
\langle \gamma, v^h - \Pi_h \Gamma(\,\cdot\,; v^h) \rangle_X = \langle \gamma, \Pi_h^\perp \Gamma(\,\cdot\,; v^h) \rangle_X  = \langle \Pi_h^\perp \gamma, \Pi_h^\perp \Gamma(\,\cdot\,; v^h) \rangle_X, \quad \forall \; v^h \in X^h.
\ee
The estimation of the RHS in the above identity will be considered in Lemma~\ref{Lem_EL_identity_RHS}; and the left-hand-side (LHS) will be considered in Lemma~\ref{Lem_EL_identity_LHS}.

\bl \label{Lem_EL_identity_RHS}
Let $\Gamma$ and $\gamma$ be defined in \eqref{Eq_Gamma} and \eqref{Eq_gamma}, respectively. Let $\overline{h}$ be as specified in Theorem~\ref{Thm_errest}. Then, there exists a constant $C > 0$ independent of $h$, such that for any $h \in (0, \overline{h})$, it holds that
\be \label{Eq_Gamma_est_goal}
|\langle \Pi_h^\perp \gamma, \Pi_h^\perp \Gamma(\,\cdot\,; v^h) \rangle_X | \le C  \|(w_*[y_*^h])''\|_{L^2(0,T; \mathbb{R}^d)} \, \|v^h\|_X h^2, \quad \forall \; v^h \in X^h.
\ee
\el
\begin{Proof}
The result follows essentially from the estimate \eqref{Eq_projection_error_linearFE_c} in Lemma~\ref{Lem_projection_linearFE}. First note that since $\gamma = y_*^h - w_*[y_*^h]$ and $y_*^h \in X^h$, we have
$\Pi_h^\perp \gamma = \Pi_h^\perp (y_*^h - w_*[y_*^h]) = -\Pi_h^\perp w_*[y_*^h]$. This together with \eqref{Eq_projection_error_linearFE_c} implies that
\be \label{Eq_gamma_est0}
\|\Pi_h^\perp \gamma\|_X  \le C h \|(w_*[y_*^h])''\|_{L^2(0, T; \mathbb{R}^d)}.
\ee
Again by \eqref{Eq_projection_error_linearFE_c}, we have also
\be \label{Eq_Gamma_est0}
\|\Pi_h^\perp \Gamma(\,\cdot\,; v^h)\|_X  \le C h \|\Gamma''(\,\cdot\,; v^h)\|_{L^2(0, T; \mathbb{R}^d)}.
\ee
It remains to estimate $\|\Gamma''(\,\cdot\,; v^h)\|_{L^2(0, T; \mathbb{R}^d)}$.

Since $\Gamma'(t; v^h) = \DF(y^h_*(t)) v^h(t)$, we get for almost every $t$ in $[0, T]$ that
\bes
\Gamma''(t; v^h) = [\DtF(y^h_*(t)) (y^h_*(t))'] v^h(t) + \DF(y^h_*(t))(v^h(t))'.
\ees
Then,
\be \label{Eq_Gamma_est1}
\|\Gamma''(\,\cdot\,;  v^h)\|^2_{L^2(0,T; \mathbb{R}^d)}  = \int_0^T \big\|[\DtF(y^h_*(t)) (y^h_*(t))'] v^h(t) + \DF(y^h_*(t))(v^h(t))' \big\|^2 \ \d t.
\ee
Note that
\bea \label{Eq_Gamma_est2a}
\int_0^T \big\|[\DtF&(y^h_*(t))(y^h_*(t))'] v^h(t) \big\|^2 \ \d t \\
& \le \int_0^T \|\DtF(y^h_*(t)) (y^h_*(t))'\|^2_{\mathrm{op}} \|v^h(t)\|^2 \ \d t \\
& \le \Big(\max_{t \in [0, T]} \|v^h(t)\|^2 \Big) \int_0^T  \|\DtF(y^h_*(t)) (y^h_*(t))'\|^2_{\mathrm{op}} \ \d t,
\eea
where $\|\cdot\|_{\mathrm{op}}$ is the operator norm of a matrix (cf.~Table \ref{tab_symbol_list}). To proceed further, note that
for any $z \in \mathbb{R}^d$, the Hessian $\DtF(z)$ is a bounded linear map from $\mathbb{R}^d$ into $L(\mathbb{R}^d, \mathbb{R}^d)$. We denote by $\vertiii{\DtF(z)}$ the operator norm of $\DtF(z)$. Namely,
\be
\vertiii{\DtF(z)} = \sup_{w \in \mathbb{R}^d, \|w\| = 1} \|\DtF(z)w \|_{\mathrm{op}}.
\ee
Since $F$ is assumed to be $C^3$, $\vertiii{\DtF(z)}$ is bounded for all $z$ on any bounded set of $\mathbb{R}^d$. Let $\mathfrak{B}$ be the subset in $\mathbb{R}^d$ defined by \eqref{Def_set_B}. We get then
\bea \label{Eq_Gamma_est2b}
\int_0^T  \|\DtF(y^h_*(t)) (y^h_*(t))'\|^2_{\mathrm{op}} \ \d t & \le \sup_{z \in \mathfrak{B}} \vertiii{\DtF(z)}^2 \int_0^T  \|(y^h_*(t))'\|^2 \ \d t \\
&\le \sup_{z \in \mathfrak{B}} \vertiii{\DtF(z)}^2 \|y^h_*\|_X^2 \\
&\le \sup_{z \in \mathfrak{B}} \vertiii{\DtF(z)}^2 (r + \|y_*\|_X)^2,
\eea
 where the last inequality follows since $y^h_* \in B(y_*,r)$ for all $h \in (0, \overline{h})$; cf.~Assumption {\bf (A2)}. Using \eqref{Eq_Gamma_est2b} in \eqref{Eq_Gamma_est2a} and noticing that $\max_{t \in [0, T]} \|v^h(t)\| \le \mathfrak{C} \|v^h\|_X$ (cf.~\eqref{Eq_embedding_const}), we get
\be \label{Eq_Gamma_est3}
\int_0^T \big\|[\DtF(y^h_*(t))(y^h_*(t))'] v^h(t) \big\|^2 \ \d t \le \Big(\mathfrak{C} \sup_{z \in \mathfrak{B}} \vertiii{\DtF(z)} (r + \|y_*\|_X)\|v^h\|_X \Big)^2.
\ee
Note also that
\bea \label{Eq_Gamma_est4}
\int_0^T \big\|\DF(y^h_*(t))(v^h(t))' \big\|^2 \ \d t & \le \int_0^T \|\DF(y^h_*(t))\|^2_{\mathrm{op}} \|(v^h(t))'\|^2 \ \d t \\
& \le \sup_{z \in \mathfrak{B}} \|\DF(z)\|^2_{\mathrm{op}} \|v^h\|_X^2.
\eea
By using \eqref{Eq_Gamma_est3} and \eqref{Eq_Gamma_est4}, we get from \eqref{Eq_Gamma_est1} that
\be \label{Eq_Gamma_est5}
\|\Gamma''(\,\cdot\,;  v^h)\|_{L^2(0,T; \mathbb{R}^d)} \le C \|v^h\|_X,
\ee
where the constant $C>0$ depends on $\sup_{z \in \mathfrak{B}} \|\DF(z)\|_{\mathrm{op}}$, $\sup_{z \in \mathfrak{B}} \vertiii{\DtF(z)}$, $\|y_*\|_X$ and the embedding constant $\mathfrak{C}$, but is independent of $h$. The desired estimate \eqref{Eq_Gamma_est_goal} follows now from \eqref{Eq_gamma_est0}, \eqref{Eq_Gamma_est0} and \eqref{Eq_Gamma_est5}.
\qed
\end{Proof}

The term $\langle \gamma, v^h - \Pi_h \Gamma(\,\cdot\,; v^h) \rangle_X$ on the LHS of Eq.~\eqref{Eq_EL_identity_rearranged} can be handled using the following lemma.
\bl \label{Lem_EL_identity_LHS}
Let $\widetilde{\mathfrak{C}} > 0$ be the embedding constant between $L^2(0,T; \mathbb{R}^d)$ and $X$. If
\be
\widetilde{\mathfrak{C}} \sup_{z \in \mathfrak{B}}\|\DF(z)\|_{\mathrm{op}} < 1,
\ee
then there exists a unique $\widehat{v}^h \in X^h$ satisfying
\be \label{Eq_vh_choice}
\widehat{v}^h - \Pi_h \Gamma(\,\cdot\,; \widehat{v}^h)  =  y_*^h - w_*^h[y_*^h].
\ee
Moreover, it holds that
\be \label{Eq_vhat_est}
\|\widehat{v}^h\|_X \le \frac{1}{1- \widetilde{\mathfrak{C}} \sup_{z \in \mathfrak{B}}\|\DF(z)\|_{\mathrm{op}}}\|y_*^h - w_*^h[y_*^h]\|_X.
\ee
\el

\begin{Proof}
Note that the map $\Psi \colon v^h \rightarrow \Pi_h \Gamma(\,\cdot\,; v^h)$ is a bounded linear map from $X^h$ onto itself.
To guarantee the existence of a unique $\widehat{v}^h$ that satisfies \eqref{Eq_vh_choice}, we only need to show that $\mathrm{Id}_{X^h} - \Psi$ is invertible. By the definition of $\Gamma$ in \eqref{Eq_Gamma}, we have
\beas
\|\Psi(v_h)\|_X = \|\Pi_h \Gamma(\,\cdot\,; v^h)\|_X & \le  \max_{t\in[0,T]}\|\DF(y^h_*(t))\|_{\mathrm{op}}\|v^h\|_{L^2(0,T; \mathbb{R}^d)}  \\
& \le \widetilde{\mathfrak{C}}   \sup_{z \in \mathfrak{B}}\|\DF(z)\|_{\mathrm{op}} \|v^h\|_X, \quad \forall \; v^h \in X^h.
\eeas
We get then
\be \label{Eq_Psi_est}
\|v_h - \Psi(v_h)\|_X \ge (1 - \widetilde{\mathfrak{C}} \sup_{z \in \mathfrak{B}}\|\DF(z)\|_{\mathrm{op}} ) \|v_h\|_X, \quad \forall \; v^h \in X^h.
\ee
Since $\widetilde{\mathfrak{C}} \sup_{z \in \mathfrak{B}}\|\DF(z)\|_{\mathrm{op}} < 1$ by our assumption, it follows that the operator norm of $\mathrm{Id}_{X^h} - \Psi$ is bounded below by $1 - \widetilde{\mathfrak{C}} \sup_{z \in \mathfrak{B}}\|\DF(z)\|_{\mathrm{op}}$. It is thus indeed invertible. The estimate \eqref{Eq_vhat_est} follows directly from \eqref{Eq_vh_choice} and \eqref{Eq_Psi_est}.
\qed
\end{Proof}

We are now in position to prove Proposition~\ref{Prop_errest_part1}.
\begin{Proof}[Proof of  Proposition~\ref{Prop_errest_part1}]
With $\widehat{v}^h$ chosen so that \eqref{Eq_vh_choice} holds and recall the definition of $\gamma$ given by \eqref{Eq_gamma}, we get
\be \label{Eq_EL_identity_LHS}
\langle \gamma, \widehat{v}^h - \Pi_h \Gamma(\,\cdot\,; \widehat{v}^h) \rangle_X = \langle y_*^h - w_*[y_*^h], y_*^h - w_*^h[y_*^h] \rangle_X.
\ee
By rewriting $y_*^h - w_*[y_*^h]$ as $y_*^h - w_*[y_*^h] = \big(y_*^h - w^h_*[y_*^h] \big) + \big(w^h_*[y_*^h] - w_*[y_*^h] \big)$, and recalling from \eqref{Eq_nonlinear_modified_FEM_soln} that $w^h_*[y_*^h] - w_*[y_*^h] = - \Pi_h^\perp w_*[y_*^h]$ lives in the orthogonal complement of $X^h$, we obtain
\bes
\langle y_*^h - w_*[y_*^h], y_*^h - w_*^h[y_*^h] \rangle_X = \langle y_*^h - w^h_*[y_*^h], y_*^h - w_*^h[y_*^h] \rangle_X = \|y_*^h - w^h_*[y_*^h]\|_X^2.
\ees
Using this identity in \eqref{Eq_EL_identity_LHS}, we get
\be \label{Eq_EL_identity_LHS_v2}
\langle \gamma, \widehat{v}^h - \Pi_h \Gamma(\,\cdot\,; \widehat{v}^h) \rangle_X = \|y_*^h - w^h_*[y_*^h]\|_X^2.
\ee
Note also that by taking $v^h = \widehat{v}^h$ in \eqref{Eq_EL_identity_rearranged}, we have
\be \label{Eq_EL_identity_rearranged2}
\langle \gamma, \widehat{v}^h - \Pi_h \Gamma(\,\cdot\,; \widehat{v}^h) \rangle_X  = \langle \Pi_h^\perp \gamma, \Pi_h^\perp \Gamma(\,\cdot\,; \widehat{v}^h) \rangle_X.
\ee
Now, it follows from \eqref{Eq_EL_identity_LHS_v2}, \eqref{Eq_EL_identity_rearranged2}, and \eqref{Eq_Gamma_est_goal}  that
\be \label{Eq_err_part1}
\|y_*^h - w^h_*[y_*^h]\|_X^2 \le C  \|(w_*[y_*^h])''\|_{L^2(0,T; \mathbb{R}^d)} \, \|\widehat{v}^h\|_X h^2.
\ee
Recall also from Lemma~\ref{Lem_EL_identity_LHS} that
\be \label{Eq_vh_est}
\|\widehat{v}^h\|_X \le \frac{1}{1 - \widetilde{\mathfrak{C}} \sup_{z \in \mathfrak{B}}\|\DF(z)\|_{\mathrm{op}}} \|y_*^h - w_*^h[y_*^h]\|_X.
\ee
The desired estimate \eqref{Eq_errest_part1_goal} follows from \eqref{Eq_err_part1} and \eqref{Eq_vh_est}.
\qed
\end{Proof}

In the proof of the main theorem given below, we require an upper bound of the term $\|(w_*[y_*^h])''\|_{L^2(0, T; \mathbb{R}^d)}$ appearing on the RHS of \eqref{Eq_errest_part1_goal}. This bound should be furthermore independent of the {\tt lsfem} solution $y_*^h$. We derive now such a bound. For this purpose, we make use of the solution to the IVP
\bea \label{Eq_nonlinear_modified2}
& w' = F(y_*(t)) + f(t),  && t \in (0, T], \\
& w(0) = g.
\eea

\bl \label{Corr_errest_part1}
Let $r$, $\overline{h}$, and $\mathfrak{B}$ be as specified in Assumption {\bf (A2)}. Let $w_*[y_*]$ be the solution of \eqref{Eq_nonlinear_modified2}. Then, for each $h \in (0, \overline{h})$, the solution $w_*[y_*^h]$ of the problem \eqref{Eq_nonlinear_modified} can be estimated as:
\be \label{Eq_errest_part1_goal2}
\|(w_*[y_*^h])''\|_{L^2(0, T; \mathbb{R}^d)}  \le  \|(w_*[y_*])''\|_{L^2(0, T; \mathbb{R}^d)} + \sup_{z \in \mathfrak{B}}\|\DF(z)\|_{\mathrm{op}} \Big( r + 2 \|y_*\|_X \Big).
\ee
\el

\begin{Proof}
Note that by using the triangle inequality,
\be \label{Eq_errest_part1_goal2_est1}
\|(w_*[y_*^h])''\|_{L^2(0, T; \mathbb{R}^d)} \le \|(w_*[y_*^h] - w_*[y_*])''\|_{L^2(0, T; \mathbb{R}^d)} + \|(w_*[y_*])''\|_{L^2(0, T; \mathbb{R}^d)},
\ee
we only need to estimate the term $\|(w_*[y_*^h] - w_*[y_*])''\|_{L^2(0, T; \mathbb{R}^d)}$.
Since $w_*[y_*^h]$ and $w_*[y_*]$ are respectively the solutions to the IVPs \eqref{Eq_nonlinear_modified} and \eqref{Eq_nonlinear_modified2}, we have
\be
(w_*[y_*^h])' - (w_*[y_*])' = F(y_*^h) - F(y_*).
\ee
Then,
\be
(w_*[y_*^h] - w_*[y_*])'' = \DF(y_*^h) (y_*^h)' - \DF(y_*)(y_*)',
\ee
which leads to
\bea  \label{Eq_errest_part1_goal2_est2}
 \|(w_* & [y_*^h] - w_*[y_*])''\|_{L^2(0, T; \mathbb{R}^d)} \\
&= \|\DF(y_*^h) (y_*^h)' - \DF(y_*)(y_*)'\|_{L^2(0, T; \mathbb{R}^d)} \\
& \le \|\DF(y_*^h) (y_*^h - y_*)'\|_{L^2(0, T; \mathbb{R}^d)}  + \|(\DF(y_*^h) - \DF(y_*)) (y_*)'\|_{L^2(0, T; \mathbb{R}^d)} \\
& \le \sup_{z \in \mathfrak{B}}\|\DF(z)\|_{\mathrm{op}} \Big( \|(y_*^h - y_*)'\|_{L^2(0, T; \mathbb{R}^d)} + 2 \|(y_*)'\|_{L^2(0, T; \mathbb{R}^d)} \Big) \\
& \le  \sup_{z \in \mathfrak{B}}\|\DF(z)\|_{\mathrm{op}} \Big( \|y_*^h - y_*\|_X + 2 \|y_*\|_X \Big)  \\
& \le \sup_{z \in \mathfrak{B}}\|\DF(z)\|_{\mathrm{op}} \Big( r + 2 \|y_*\|_X \Big).
\eea
In deriving \eqref{Eq_errest_part1_goal2_est2}, we have used the facts that $\|y_*^h - y_*\|_X< r$ since $y_*^h \in B(y_*,r)$ and that $\mathfrak{B}$ contains both $y_*(t)$ and $y^h_*(t)$ for all $t \in [0, T]$; see \eqref{Def_set_B} and \eqref{Eq_LSFEM_in_B}. The desired result \eqref{Eq_errest_part1_goal2} follows from \eqref{Eq_errest_part1_goal2_est2} and \eqref{Eq_errest_part1_goal2_est1}.
\qed
\end{Proof}

We are now in position to prove the main theorem of this section.

\begin{Proof}[{\bf Proof of Theorem~\ref{Thm_errest}}]
First note that by the triangle inequality, we get
\bea \label{Eq_errest1}
 \|y_* - y_*^h\|_{L^2(0, T; \mathbb{R}^d)} & \le  \|y_* - w_*[y_*^h]\|_{L^2(0, T; \mathbb{R}^d)} \\
 & \quad +  \|w_*[y_*^h] - w_*^h[y_*^h]\|_{L^2(0, T; \mathbb{R}^d)} +  \|w_*^h[y_*^h] - y_*^h\|_{L^2(0, T; \mathbb{R}^d)}.
\eea
To estimate the first term $\|y_* - w_*[y_*^h]\|_{L^2(0, T; \mathbb{R}^d)}$ on the RHS of \eqref{Eq_errest1}, we integrate Eq.~\eqref{Eq_nonlinear} and Eq.~\eqref{Eq_nonlinear_modified} to obtain
\bea
\|y_*(t) - w_*[y_*^h](t)\| & \le \int_0^t \|F(y_*(s)) - F(y_*^h(s))\| \ \d s \\
& \le \mathrm{Lip}(F\vert_{\mathfrak{B}}) \int_0^t \|y_*(s) - y_*^h(s) \| \ \d s \\
& \le \mathrm{Lip}(F\vert_{\mathfrak{B}}) \sqrt{t} \left(\int_0^t \|y_*(s) - y_*^h(s) \|^2 \ \d s\right)^{1/2}, \qquad t \in [0, T],
\eea
where we applied H\"older's inequality in the last step above. We get in turn that
\be \label{Eq_errest2}
 \|y_* - w_*[y_*^h]\|_{L^2(0, T; \mathbb{R}^d)} \le \frac{\mathrm{Lip}(F\vert_{\mathfrak{B}}) T}{\sqrt{2}} \|y_* - y_*^h\|_{L^2(0, T; \mathbb{R}^d)}.
\ee
The second term $\|w_*[y_*^h] - w_*^h[y_*^h]\|_{L^2(0, T; \mathbb{R}^d)}$  in \eqref{Eq_errest1} can be estimated by using \eqref{Eq_nonlinear_modified_errest}, and
the last term $\|w_*^h[y_*^h] - y_*^h\|_{L^2(0, T; \mathbb{R}^d)}$ can be estimated by using \eqref{Eq_errest_part1_goal}
together with
\be \label{Eq_w_embedding_est}
\|w_*^h[y_*^h] - y_*^h\|_{L^2(0, T; \mathbb{R}^d)} \le \widetilde{\mathfrak{C}} \|w_*^h[y_*^h] - y_*^h\|_X,
\ee
where $\widetilde{\mathfrak{C}}$ denotes again the embedding constant between $L^2(0,T; \mathbb{R}^d)$ and $X$.

Gathering the above estimates for the three terms on the RHS of \eqref{Eq_errest1}, we get
\bea
 \|y_* - y_*^h\|_{L^2(0, T; \mathbb{R}^d)} &  \le  \frac{\mathrm{Lip}(F\vert_\mathfrak{B}) T}{\sqrt{2}}  \|y_* - y_*^h\|_{L^2(0, T; \mathbb{R}^d)} \\
 & \hspace{-1.5em}+ C \Big(1 + \frac{1}{1 - \widetilde{\mathfrak{C}} \sup_{z \in \mathfrak{B}}\|\DF(z)\|_{\mathrm{op}}}\Big)\|(w_*[y_*^h] )''\|_{L^2(0, T; \mathbb{R}^d)} h^2,
 \eea
where we have absorbed the factor $\widetilde{\mathfrak{C}}$ on the RHS of \eqref{Eq_w_embedding_est} into the constant $C$ when applying the estimate \eqref{Eq_errest_part1_goal}.

In the above inequality, by using the estimate \eqref{Eq_errest_part1_goal2} for $\|(w_*[y_*^h] )''\|_{L^2(0, T; \mathbb{R}^d)}$, we get after rearranging terms that
\be
\Big(1  - \frac{\mathrm{Lip}(F\vert_\mathfrak{B}) T}{\sqrt{2}} \Big)  \|y_* - y_*^h\|_{L^2(0, T; \mathbb{R}^d)} \le \widetilde{C} h^2,
\ee
where
\bea
\widetilde{C} := C \Big(1 + \frac{1}{1 - \widetilde{\mathfrak{C}} \sup_{z \in \mathfrak{B}}\|\DF(z)\|_{\mathrm{op}}}\Big) \Big( & \|(w_*[y_*])''\|_{L^2(0, T; \mathbb{R}^d)}  \\
 & + \sup_{z \in \mathfrak{B}}\|\DF(z)\|_{\mathrm{op}} \big( r + 2 \|y_*\|_X \big) \Big),
\eea
with $w_*[y_*]$ denoting the solution of \eqref{Eq_nonlinear_modified2}. Since it is assumed that $\frac{\mathrm{Lip}(F\vert_\mathfrak{B}) T}{\sqrt{2}} < 1$ (see \eqref{Eq1_in_A2}), the desired result
\eqref{Eq_errest_goal} follows by taking $C$ therein to be $ \widetilde{C} / (1  - \mathrm{Lip}(F\vert_\mathfrak{B}) T/\sqrt{2})$.
\qed
 \end{Proof}

\section{Numerics} \label{Sect_numerics}

In this section we focus on numerical aspects of {\tt lsfem} by discussing algorithmic details and confirming the analytical insight gained in previous Section~\ref{sec:theory} and \ref{Sect_theory_err_estimates} through numerical experiments. We discuss modifications such as adaptive time stepping and constrained systems, i.e., differential algebraic equations.  Our proposed method is summarized in Algorithm~\ref{alg:lsfem}.

\begin{algorithm}[H]
\caption{({\tt lsfem} for IVPs)}\label{alg:lsfem}
\begin{algorithmic}[1]
\Function{$y^h_*$ = {\tt lsfem}}{$G(\mdot,\mdot)$, $g$, $[t_0,T]$, $X^h$}
\State construct finite element basis for $y^h(\mdot;\mdot)$ and $(y^{h})' (\mdot;\mdot)$ of $X^h$\label{ln:fem}
\State choose discretization $\underline{t}$ of time interval \label{ln:disc}
\State compute
$$x_* \in \argmin_{x} \ \calJ(x)= \thf||(y^{h})'(\underline{t}; x)-G(\underline{t},y^{h}(\underline{t}; x))||^2 + \thf||y^{h}(0; x) - g||^2$$
\vspace*{-6ex}\label{ln:opt}
\State set $y^h_* = y^h(\mdot; x_*)$
\EndFunction
\end{algorithmic}
\end{algorithm}
\vspace*{-2ex}
A Matlab implementation of Algorithm~\ref{alg:lsfem} is available at {\tt \href{https://github.com/matthiaschung/lsfem}{github.com/matthiaschung/lsfem}} and is intended for reproducibility and to develop an understanding of the performance of {\tt lsfem} method for ODEs.
As inputs {\tt lsfem} requires the RHS of the first order ODE $G$, the interval of interest $[t_0,T]$, and the initial condition $y(t_0) = g$,
\bea\label{eq:ivp}
& y' = G(t,y), && t \in (t_0, T],  \\
& y(t_0) = g.
\eea
Additionally, one may select a desired finite element space $X^h$. Algorithm~\ref{alg:lsfem} return the function $y_h^* = y_h(\mdot;x_*)$ determined by the optimized finite element coefficients $x_*$ with respect to the corresponding finite element basis.

Various numerical choices need to be made in Algorithm~\ref{alg:lsfem}. First, the approximation quality of our method depends greatly on the choice of the finite element space $X^h$ and its corresponding control points $\tau$ in line~\ref{ln:fem}. Section~\ref{sec:theory} and~\ref{Sect_theory_err_estimates} provides convergence results for piecewise linear basis functions, however, we may choose higher order basis functions.  Common choices for the basis function include piecewise polynomials and polynomial splines \cite{deBoor2001,Schumaker2015}. Other options are exponential splines, which may better capture the exponential behavior manifested by certain differential equations; see \cite{Rentrop1980}. An interesting alternative are Hermite splines \cite{kreyszig11}, which are able to take advantage of derivative information provided naturally by the differential equation and reducing computational costs. Note that, the choice of the finite element basis may depend on the imposed smoothness of the underlying dynamical system, i.e., $G$. Equidistant control points may be selected if no further information on $y$ are available, however, one may also select control points if knowledge on $y$ (or its derivatives) are available.

To numerically evaluate and minimize $J$ of equation~\eqref{Eq_opt_FEM}, quadrature is required to approximate the $L^2$-norm. Hence, with respect to the quadrature rule we discretize the interval $[t_0,T]$ with $\underline{t} =[t_0,t_1,\ldots,t_{n-1},T]\t$, see line~\ref{ln:disc}. Choosing a quadrature rule (such as Gauss-Legendre and Gauss-Lobatto, \cite{Golub1967}) which is consistent with the finite element space $X^h$ may provide computational advantages. The added benefit of using such a quadrature rule is that the resulting $\ell^2$-norm approximation $\calJ$ has the potential to be exact in certain polynomial settings.

The main computational effort lies in line~\ref{ln:opt}. Line~\ref{ln:opt} defines a common (regularized) nonlinear least squares problem. Notice that, if $G$ is sufficiently smooth, gradient and Hessian based optimization methods can be utilized.

Hence, gradient based methods and also Newton type methods are natural choices (assuming sufficient smoothness of the system). However, $\nabla_{x} J$ and $\nabla_{x}^2 J$ need to be readily available or be obtained by algorithmic differentiation techniques~\cite{griewank2008evaluating}. It is worth mentioning, that {\tt lsfem} seeks for a global minimizer $y_h^*$ of~\eqref{Eq_opt_FEM}. However, for non-convex problems the proposed optimization methods may not ensure convergence to the global minimizer. Strategies to prevent local minimizer are required, e.g., multi-start or global optimization methods~\cite{horst2013handbook}.

\begin{figure}[h]
  \begin{center}
\renewcommand\figurescale{0.52}
%
%
\definecolor{mycolor1}{rgb}{0.00000,0.44700,0.74100}%
\definecolor{mycolor2}{rgb}{0.85000,0.32500,0.09800}%
\definecolor{mycolor3}{rgb}{0.92900,0.69400,0.12500}%
\definecolor{mycolor4}{rgb}{0.49400,0.18400,0.55600}%
\definecolor{mycolor5}{rgb}{0.46600,0.67400,0.18800}%
\begin{tikzpicture}

\begin{axis}[%
width=6.458in,
height=5.641in,
at={(1.083in,0.783in)},
scale only axis,
xmode=log,
xmin=0.05,
xmax=1,
xminorticks=true,
xlabel style={font=\color{white!15!black}},
xlabel={meshsize $h$ [${\rm log}_{10}$]},
ymode=log,
ymin=2.33354969513912e-13,
ymax=0.0326304697640342,
yminorticks=true,
ylabel style={font=\color{white!15!black}},
ylabel={$||y^h_* - y||_{L^2(0,T;\mathbb{R})}$ [$\log_{10}$]},
axis background/.style={fill=white},
axis x line*=bottom,
axis y line*=left,
legend style={at={(0.97,0.03)}, anchor=south east, legend cell align=left, align=left, draw=white!15!black},
scale=\figurescale, x axis line style = {thick}, y axis line style = {thick}
]
\addplot [color=mycolor1, line width=1.0pt, mark size=0.8pt, mark=*, mark options={solid, mycolor1}]
  table[row sep=crcr]{%
1	0.0326304697640342\\
0.5	0.00854019841090216\\
0.333333333333333	0.00382915082321112\\
0.25	0.00216058447927277\\
0.2	0.00138476532778239\\
0.166666666666667	0.000962395667948159\\
0.142857142857143	0.000707400298053681\\
0.125	0.000541769519305984\\
0.111111111111111	0.000428154872020254\\
0.1	0.000346857651517752\\
0.0909090909090909	0.000286691147984877\\
0.0833333333333333	0.000240920612993045\\
0.0769230769230769	0.00020529500701325\\
0.0714285714285714	0.00017702383580723\\
0.0666666666666667	0.000154213941072228\\
0.0625	0.000135544280393299\\
0.0588235294117647	0.000120070344557227\\
0.0555555555555556	0.000107102353235943\\
0.0526315789473684	9.61270563342281e-05\\
0.05	8.67561729917463e-05\\
};
\addlegendentry{$k = 1$}

\addplot [color=mycolor2, line width=1.0pt, mark size=0.8pt, mark=*, mark options={solid, mycolor2}]
  table[row sep=crcr]{%
1	0.00338562544623212\\
0.5	0.000489616211886746\\
0.333333333333333	0.000142271927624712\\
0.25	5.96867849065177e-05\\
0.2	3.04535995918953e-05\\
0.166666666666667	1.7589861559815e-05\\
0.142857142857143	1.10634909459015e-05\\
0.125	7.40567566699004e-06\\
0.111111111111111	5.19829648521309e-06\\
0.1	3.78800091698341e-06\\
0.0909090909090909	2.84510600390396e-06\\
0.0833333333333333	2.19093893637402e-06\\
0.0769230769230769	1.72291402369812e-06\\
0.0714285714285714	1.37925769523602e-06\\
0.0666666666666667	1.12125330585513e-06\\
0.0625	9.23793313953548e-07\\
0.0588235294117647	7.70109417639416e-07\\
0.0555555555555556	6.48711777169826e-07\\
0.0526315789473684	5.51547749747552e-07\\
0.05	4.72859599909054e-07\\
};
\addlegendentry{$k = 2$}

\addplot [color=mycolor3, line width=1.0pt, mark size=0.8pt, mark=*, mark options={solid, mycolor3}]
  table[row sep=crcr]{%
1	0.000202183391105212\\
0.5	2.32904464670047e-05\\
0.333333333333333	6.13925571326068e-06\\
0.25	1.98053469371163e-06\\
0.2	8.53606226589755e-07\\
0.166666666666667	4.17652307472761e-07\\
0.142857142857143	2.29487981236018e-07\\
0.125	1.35838535734138e-07\\
0.111111111111111	8.55637985267875e-08\\
0.1	5.65006508146696e-08\\
0.0909090909090909	3.88028978931251e-08\\
0.0833333333333333	2.75186658164752e-08\\
0.0769230769230769	2.00544994470919e-08\\
0.0714285714285714	1.49573475416082e-08\\
0.0666666666666667	1.13814219621748e-08\\
0.0625	8.8129295115025e-09\\
0.0588235294117647	6.92975385672289e-09\\
0.0555555555555556	5.52372173313448e-09\\
0.0526315789473684	4.45685236870302e-09\\
0.05	3.63554315284495e-09\\
};
\addlegendentry{$k = 3$}

\addplot [color=mycolor4, line width=1.0pt, mark size=0.8pt, mark=*, mark options={solid, mycolor4}]
  table[row sep=crcr]{%
1	9.76821015938524e-06\\
0.5	9.44614402472246e-07\\
0.333333333333333	3.47146887969513e-07\\
0.25	8.20255205439677e-08\\
0.2	2.89194787914358e-08\\
0.166666666666667	1.15664390181862e-08\\
0.142857142857143	5.46097611519536e-09\\
0.125	2.80196850460666e-09\\
0.111111111111111	1.5667573352963e-09\\
0.1	9.26734778355406e-10\\
0.0909090909090909	5.77482950729021e-10\\
0.0833333333333333	3.74385439626387e-10\\
0.0769230769230769	2.51413163418983e-10\\
0.0714285714285714	1.737987138036e-10\\
0.0666666666666667	1.23255855037096e-10\\
0.0625	8.93541254695691e-11\\
0.0588235294117647	6.60518769676769e-11\\
0.0555555555555556	4.96727170608683e-11\\
0.0526315789473684	3.79340228905776e-11\\
0.05	2.93713873207182e-11\\
};
\addlegendentry{$k = 4$}

\addplot [color=mycolor5, line width=1.0pt, mark size=0.8pt, mark=*, mark options={solid, mycolor5}]
  table[row sep=crcr]{%
1	3.98268938643709e-07\\
0.5	3.30897017781375e-08\\
0.333333333333333	1.75577249373246e-08\\
0.25	3.35878162438995e-09\\
0.2	9.15004727585026e-10\\
0.166666666666667	3.08935486517417e-10\\
0.142857142857143	1.24019677601377e-10\\
0.125	5.58001694404935e-11\\
0.111111111111111	2.7681455587558e-11\\
0.1	1.47364975279712e-11\\
0.0909090909090909	8.34425374731115e-12\\
0.0833333333333333	4.95767115029332e-12\\
0.0769230769230769	3.0728046123035e-12\\
0.0714285714285714	1.97214466080051e-12\\
0.0666666666666667	1.30540765546396e-12\\
0.0625	8.87089425094118e-13\\
0.0588235294117647	6.17174521389902e-13\\
0.0555555555555556	4.38453598567911e-13\\
0.0526315789473684	3.171383867142e-13\\
0.05	2.33354969513912e-13\\
};
\addlegendentry{$k = 5$}

\end{axis}
\end{tikzpicture}%
  \end{center}
  \caption{Illustration of convergence rates of {\tt lsfem} methods for B-spline bases of degree $k$ ($k = 1,\ldots, 5$) in comparison to the mesh sizes $h$ for the ODE $y' = -y$ with $t \in [0,1]$ and $y(0) = 1$. The numerically observed convergence rates are 1.9972 (for $k=1$),  3.0066 (for $k=2$), 3.9204 (for $k=3$), 4.9456 (for $k=4$) and 5.9409 (for $k=5$).}\label{fig:converenceRates}
\end{figure}
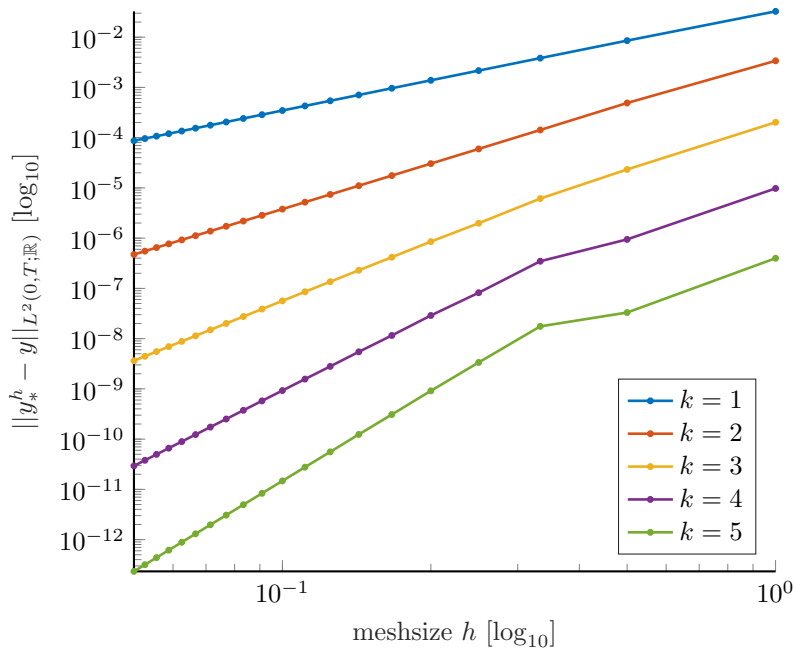

\medskip
\paragraph{\it Rate of convergence for higher order finite elements}
To illustrate and empirically valid the convergence rates discussed in Section~\ref{sec:theory} and~\ref{Sect_theory_err_estimates}, we first consider the linear initial value problem $y' = -y$, with $t \in [0,1]$ and $y(0) = 1$. We use a B-spline bases for $X^h$ of degree $k = 1,\ldots, 5$ with varying equidistant discretization of the finite elements, i.e., $h = 1/N$ with $N = 1,\ldots, 20$ to compute the finite element approximation $y_*^h$. Figure~\ref{fig:converenceRates} depicts the errors  $\| y_*^h - y\|$ with respect to the varying mesh sizes $h$ in a log-log space. The slopes of each graph reveal the power of the expected convergence rates of our method. For instance the slope using linear B-splines is 1.9972 confirming the quadratic convergence rate (Theorem~\ref{Thm_errest}). The other rates are 3.0066 (for $k=2$), 3.9204 (for $k=3$), 4.9456 (for $k=4$) and 5.9409 (for $k=5$), respectively. These results lead us to conjecture that the optimal {\tt lsfem} error bounds scale like  $\calO(h^{k+1})$ for finite element bases of degree $k$. By inspecting the proofs in Section~\ref{Sect_theory_err_estimates}, we expect many of the ingredients presented there to be extended naturally, although some aspects such as higher-order analogues of Lemma~\ref{Corr_errest_part1} may require additional efforts. We plan to address such an extension in a future work.

\begin{figure}[h]
  \begin{center}
      \renewcommand\figurescale{0.60}%
    \input{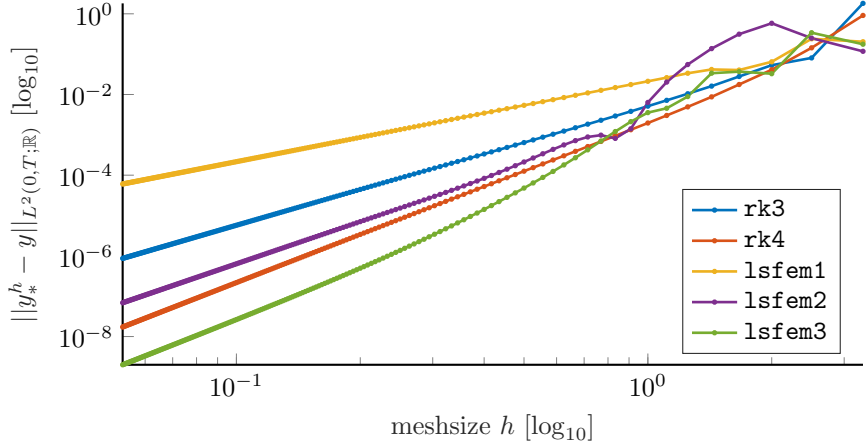}
  \end{center}
  \caption{Illustration of convergence rates of {\tt lsfem} utilizing B-spline bases of degree 1, 2, and 3 in comparison to Runge-Kutta 3 and 4 for the logistic equation $y' = y(1-y)$, with $y(0) = 1/10$ and $t \in [0, 10]$. The numerically observed convergence rates are 2.001, 3.4928, and 4.0470 for the {\tt lsfem} methods and 2.9731, 3.9820 for the Runge-Kutta methods, respectively.}\label{fig:converenceRatesLogistic}
\end{figure}

To confirm the results for nonlinear ODEs we consider the simple logistic growth $y' = y(1-y)$ with $y(0) = 1/10$ and $t \in [0,10]$. We compare our {\tt lsfem} method with B-spline bases of degree 1, 2, and 3 ({\tt lsfem1}, {\tt lsfem2}, and {\tt lsfem3}) to Runge-Kutta 3 ({\tt rk3}) and Runge-Kutta 4 ({\tt rk4}), see Figure~\ref{fig:converenceRatesLogistic}. The numerically observed convergence rates for this logistic growth model are 2.001, 3.4928, and 4.0470 for the {\tt lsfem} methods and 2.9731, 3.9820 for the Runge-Kutta methods, respectively. The observed rates for {\tt lsfem} confirm again the obtained theoretical estimate for degree $k = 1$ case and corroborate the conjectured optimal bound $\calO(h^{k+1})$ for higher-degree bases (with {\tt lsfem2} providing actually better rate than conjectured for this particular example). One can also compare {\tt lsfem3} with {\tt rk4}, since both methods show a convergence rate close to the theoretical rate $r=4$. Figure~\ref{fig:converenceRatesLogistic} reveals that the constant $C$ in the associated error bound $C h^4$ is smaller in the case of {\tt lsfem3} than that of {\tt rk4} for the considered example.

\medskip
\paragraph{\it Linear ODEs}
In case of linear ODEs the {\tt lsfem}'s main computational burden of solving the optimization problem in line~\ref{ln:opt} of Algorithm~\ref{alg:lsfem} simplifies to a linear least-squares problem whose solution can be obtained e.g., by solving the associated linear normal equations.

More precisely, let us consider the $n$ dimensional initial value problem
\begin{equation}\label{eq:nonhonmo}
	y'(t) = A(t)y(t) + b(t) \quad \mbox{and } y(t_0) = g \mbox{ and } t\in(t_0,T].
\end{equation}
Assuming we choose the same finite element basis for each state
\begin{equation}
		\phi(t)=[\phi_1(t),\phi_2(t),...,\phi_m(t)]\t,
\end{equation}
then the function of the finite element space are given by $y^h(t, x)=(\phi(t)\t\kron I_n)x$ with some coefficients
\begin{equation}
		x=\left[x_1^1,\ldots,x_1^n, x_2^1,\ldots,x_2^n,\ldots,x_m^1,\ldots,x_m^n\right]\t
\end{equation}
where $\kron$ denotes the Kronecker product and $I_n$ the identity matrix. The least-squares problem now reads
\begin{equation}
	\min_x \ \thf\int_{t_0}^T r(t,x)\t r(t,x) \ \d t + \thf \norm[2]{y^h(t_0,x) - g}^2
\end{equation}
with $r(t,x)=Z(t)x-b(t)$, where $Z(t)=\phi'(t)\t\kron I_n- A(t)(\phi(t)\t\kron I_n)$. With the further abbreviations
\begin{equation}
	 Q = \int_{t_0}^T  Z(t)\t Z(t)\ \d t, \quad
   p = \int_{t_0}^T  Z(t)\t b(t)\ \d t, \quad \mbox{and} \quad
   R = \phi(t_0)\t\kron I_n.
\end{equation}
The {\tt lsfem} solution of \eqref{eq:nonhonmo} is obtained by the normal equations
\begin{equation}
  (Q + R\t R) x_* = p + R\t g,
\end{equation}
and $y^h_*(t) = (\phi(t)\t\kron I_n)(Q + R\t R)^{-1}(p + R\t g)$, assuming  $Q + R\t R$ is invertible. Hence standard linear algebra libraries may be utilized to solve a linear system of differential equations efficiently.

\medskip
\paragraph{\it Adaptive discretization of the finite elements}
So far we have not discussed how to select control points $\tau$ of our finite element space $X^h$ and assumed they are pre-selected, e.g., equidistant. Alternatively, control points may be selected adaptively by (for simplicity) repeated evaluation of line~\ref{ln:opt} with refined control points $\tau$. Notice, {\tt lsfem} naturally provides error estimates through the residuals $r =(y^h)'(\underline{t}; x_*)-G(\underline{t},y^h(\underline{t}; x_*))$. In its simplest form new control points $\tau_i$'s may be introduced by selecting discretization points $\underline{t}_j$'s at locations with large residuals $r_j$.

\begin{figure}
    \centering
    \renewcommand\figurescale{0.47}%
    \input{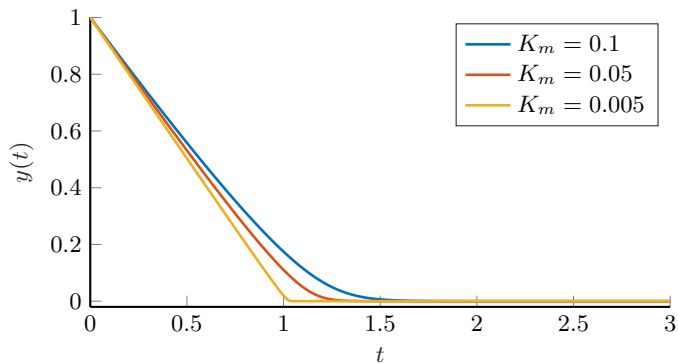}
     \caption{Solution of Eq.~\eqref{eq:kinetics} for three different values of $K_m$ and with $y(0) = 1$. By decreasing $K_m$ the rapid change of dynamics of $y$ become more prominent around $t = 1$.}
    \label{fig:kinetic}
\end{figure}

\begin{figure}
    \centering
    \renewcommand\figurescale{0.57}%
    \input{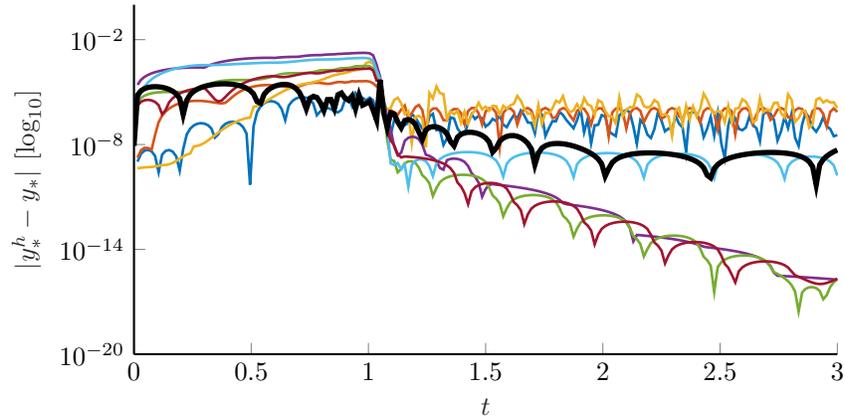}
    \caption{Absolute error between various ODE methods and the true solution for Eq.~\eqref{eq:kinetics} with $K_m = 0.005$ and $y(0) = 1$. Our {\tt lsfem} method is highlighted in black bold, other standard Matlab ODE solvers in default settings are illustrated as labelled in Figure~\ref{fig:kineticDiscretization}.}
    \label{fig:kinetcError}
\end{figure}

\begin{figure}
    \centering
    \renewcommand\figurescale{0.57}%
    \input{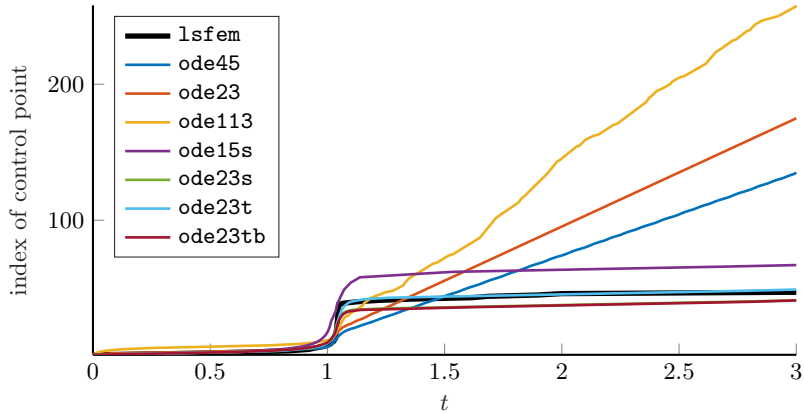}
    \caption{Time location of the corresponding control point for {\tt lsfem} (black bold) in comparison with time location of the discretization index of FDMs, for Eq.~\eqref{eq:kinetics} with $K_m = 0.005$ and $y(0) = 1$. Notice that all methods, especially {\tt lsfem} and all the stiff FDM solvers, select small step sizes around $t=1$ where the kink in the true solution occurs (cf.~Figure~\ref{fig:kinetic}). At the same time, the three non-stiff solvers {\tt ode45}, {\tt ode23}, and {\tt ode113} also require a fine discretization beyond $t=1$ even though the solution dynamics is ``quiescent'' there (see again Figure~\ref{fig:kinetic}). Note also that the curve for {\tt ode23s} almost overlaps with {\tt ode23tb}. All Matlab ODE solvers are in default settings with default error controls.}
    \label{fig:kineticDiscretization}
\end{figure}

We illustrate this procedure by considering the model
\begin{equation}\label{eq:kinetics}
  y' = -\frac{y}{K_m+y},
\end{equation}
with $y(0) = 1$, $t\in [0,3]$. The solution is implicitly given by $y + K_m \ln y = 1-t$. Equation~\eqref{eq:kinetics} is a simplified model frequently appearing in the field of enzyme kinetics, \cite{segel1988validity,eilertsen2020quasi}. Here $K_m>0$ refers to the Michaelis-Menten constant determining the reaction rate. For small $K_m$, e.g., $K_m = 0.005$, stiff ODE~\eqref{eq:kinetics} exhibits a decay with a sharp ``kink'' around $t=1$, see Figure~\ref{fig:kinetic}, resulting in difficulties for numerical ODE solvers to preserve non-negative concentrations. Adaptive refinement of this particular area is crucial. We initialize {\tt lsfem} with four equidistant control points (order 3 and 8 Gauss Legendre points) and refine the control points $\tau$ of our finite element basis until each residual element reaches an absolute tolerance of $10^{-4}$. The error of {\tt lsfem} with respect to the true solution and in comparison of standard Matlab is depicted in Figure~\ref{fig:kinetcError} while Figure~\ref{fig:kineticDiscretization} show the number of discretization points vs.~the location of these points. We observe that {\tt lsfem} adaptively adds control points around 1 and maintains an absolute error below $5\cdot 10^{-5}$ throughout the time interval while requiring 47 control points.

In comparison, stiff ODE solvers such as {\tt ode15s}, {\tt ode23s}, {\tt ode23t}, and {\tt ode23tb} also adaptively refine around $t=1$ and the total time steps used are similar to the number of control points for {\tt lsfem} except {\tt ode15s}, with the latter taking more time steps than the other stiff solvers as shown in Figure~\ref{fig:kineticDiscretization}. Note also that the curve for {\tt ode23s} almost overlaps with {\tt ode23tb} in  Figure~\ref{fig:kineticDiscretization} and is thus not visible. In terms of errors, the {\tt lsfem} performs better over the interval $[0,1]$ in which the dynamics is ``non-trival'',  and the stiff solvers perform better over the interval $[1, 3]$ in which the dynamics is ``quiescent''; see Figure~\ref{fig:kinetic} and Figure~\ref{fig:kinetcError}. Overall, {\tt lsfem} performs similar to stiff FDM solvers, while maintaining a slightly lower $L^2$-error. Meanwhile, non-stiff ODE solvers such as {\tt ode23}, {\tt ode45}, and {\tt ode113} refine less around $t=1$ but need a significantly more number of discretization points beyond $t>1$ to maintain numerical accuracy.

\section{Conclusion and discussion} \label{Sect_conclusion}

In this work, we considered the least-squares finite element method (\texttt{lsfem}) for systems of nonlinear ordinary differential equations and established under suitable conditions an optimal error estimate for this method when piecewise linear elements are used (Theorem~\ref{Thm_errest}). In contrast to the ``localization'' nature of finite difference methods, the {\tt lsfem} aims to find an optimal approximate solution within a given subspace that minimizes an objective function over the whole time interval of integration. The {\tt lsfem} can thus be less prone to the accumulation of local discretization errors compared to finite difference methods.

As reviewed in Section \ref{sec:intro}, a key ingredient in our derivation of the optimal estimate is a geometric (orthogonality) property derived from the first-order optimality condition associated with the minimizers of the underlying optimization problems; see Eq.~\eqref{Eq_EL_identity_intro}. Numerical results presented in Section~\ref{Sect_numerics} not only support our main theoretical result presented in Theorem~\ref{Thm_errest}, but also provide strong indication that error bound of the form $\calO(h^{k+1})$ will hold if higher-order spline basis elements of degree $k$ ($k\ge2$) were used. In Section~\ref{Sect_numerics}, we also discussed details related to the associated algorithmic aspects (Algorithm~\ref{alg:lsfem}) as well as suitable modifications for adaptive mesh refinement to handle ODEs whose solutions may experience abrupt local changes. It is also worth mentioning that residual neural networks appear to be of such nature and we will dedicate future research towards such applications, \cite{haber2017stable,gunther2021spline}.

Finally, we mention that the procedure presented in Algorithm~\ref{alg:lsfem} can be easily adapted to handle a broad class of differential algebraic equations (DAEs) \cite{Ascher1998,Hairer1996,reis2014surveys} as well. One just needs to add the corresponding algebraic equations as constraints to the associated optimization problems. The numerical setup can also be easily extended to handle ODE boundary value problems and delay differential equations. We plan to address these extensions in future communications.

\section*{Acknowledgments}
The authors wish to thank Micka\"el D. Chekroun and Tao Lin for stimulating and fruitful discussions on various topics related to this work. The authors are also grateful to Jeff Borggaard for suggesting the example presented in the Introduction section. This work has been partially supported by the National Science Foundation grants DMS-1723005 (M Chung) and DMS-2108856 (H Liu).

\appendix

\section{Proof Theorem \ref{Thm_convergence}} \label{Sect_conv_thm_proof}

We start by rewriting the respective first-order optimality condition associated with the minimization problems \eqref{Eq_opt} and \eqref{Eq_opt_FEM} into an abstract form $u+\mathcal{T} \circ \mathcal{G}(u) = 0$ for \eqref{Eq_opt} and $u^h+ \Pi_h \mathcal{T} \circ \mathcal{G}(u^h) = 0$ for \eqref{Eq_opt_FEM}, where $\mathcal{T}$ is a bounded linear operator and $\mathcal{G}$ is a smooth nonlinear operator defined below, and $\Pi_h$ is the orthogonal projection onto $X^h$ appearing in \eqref{Eq_Pih_conv}. These equations are in the same functional forms dealt with in \cite[Theorem 3.3, p.307]{girault2012finite}. Once this reformulation is done, we just need to check that all conditions required in \cite[Theorem 3.3, p.307]{girault2012finite} are satisfied and that the solutions for $u^h+ \Pi_h \mathcal{T} \circ \mathcal{G}(u^h) = 0$ guaranteed by \cite[Theorem 3.3, p.307]{girault2012finite} are indeed minimizers of \eqref{Eq_opt_FEM}. The reformulation is dealt with in Lemma~\ref{Lem:operator_form}, and the remaining steps are carried out afterwards.

As a preparation, we define the linear operator $\mathcal{T}$ to be the solution operator associated with the special case of \eqref{Eq_nonlinear} in which $F$ is identically zero. That is,
\be \label{Eq_def_T}
\mathcal{T}  \colon L^2(0, T; \mathbb{R}^d)  \times \mathbb{R}^d  \rightarrow X, \; \mathcal{T}(f,g)  = \ell, \text{ with }  \ell(t) = g + \int_{0}^t \! f(s) \, \d s, \; \forall\, t \in [0, T].
\ee
Note that $\mathcal{T}$ thus defined is obviously a bounded linear operator. It follows directly from the Riesz representation theorem (see e.g.,~\cite[Appendix D]{Evans10} and \cite[Section 16.2]{Royden2010real}), that $\ell=\mathcal{T}(f,g)$ defined above is the unique element in $X$ satisfying
\be \label{Eq_T_equiv_formulation}
\langle \ell,v \rangle_X = \int_0^T \langle f(t) ,  v'(t) \rangle \ \d t + \langle g, v(0) \rangle,  \qquad \forall \; v \in X.
\ee
Note also that this identity is the first-order necessary condition for $\ell$ to be a solution of \eqref{Eq_opt} when $F$ is absent; cf.~\eqref{Eq_1st_order_cond} below.

Assume that $F$ is $C^1$ smooth. For any $(f,g)$ in $L^2(0, T; \mathbb{R}^d) \times \mathbb{R}^d$, we define now a nonlinear operator $\mathcal{G}$ as follows:
\be \label{Eq_def_G}
\mathcal{G} \colon X \rightarrow L^2(0, T;  \mathbb{R}^d) \times \mathbb{R}^d, \qquad \mathcal{G}(u) = (\widetilde{f}(u), \widetilde{g}(u)),
\ee
where for each $u$ in $X$, $(\widetilde{f}(u), \widetilde{g}(u))$ is defined by
\bea \label{Eq_f_g_tilde}
& [\widetilde{f}(u)](t) = -F(u(t)) - f(t) - \!\int_t^T[\DF(u(s))]\t \big(u'-F(u) - f(s)\big)\, \d s,  \; \forall \,  t \in [0, T], \\
& \widetilde{g}(u) = -g - \int_0^T [\DF(u(s))]\t \big(u'-F(u) - f(s)\big)\, \d s.
\eea
To see that $(\widetilde{f}(u), \widetilde{g}(u))$ thus defined is indeed an element in $L^2(0, T;  \mathbb{R}^d) \times \mathbb{R}^d$, note that since $u\in X$ and $F$ is assumed to be $C^1$, we have $w=u' -  F(u) - f \in L^2(0, T; \mathbb{R}^d)$ and $\DF(u)$ is continuous on $[0,T]$. One can then show that $\int_{0}^T (\DF(u))\t w\, \d s$ is finite and the function $t \mapsto \psi(t):= \int_{t}^T (\DF(u))\t w\, \d s$ is in $L^2(0, T; \mathbb{R}^d)$. As a result, $\widetilde{f}$ maps $X$ into $L^2(0, T;  \mathbb{R}^d)$ and $\widetilde{g}$ maps $X$ into $\mathbb{R}^d$.

The rationale behind the definition of $\mathcal{G}$ will become apparent in the proof of the following lemma.

\bl \label{Lem:operator_form}
Consider the IVP \eqref{Eq_nonlinear}. Assume that $f \in L^2(0, T; \mathbb{R}^d)$ and $F: \mathbb{R}^d \rightarrow \mathbb{R}^d$ is~$C^1$. Let $\mathcal{T}$ and $\mathcal{G}$ be defined by \eqref{Eq_def_T} and \eqref{Eq_def_G}, respectively. Then,  any strong solution of \eqref{Eq_nonlinear} also satisfies the following nonlinear problem
\be \label{Eq_1st_order_cond_operator_form}
u + \mathcal{T} \circ \mathcal{G}(u) = 0.
\ee
Similarly, denoting by $\Pi_h \colon X \rightarrow X^h$ the orthogonal projection onto $X_h$, any solution of \eqref{Eq_opt_FEM} also satisfies
\be \label{Eq_1st_order_cond_operator_form_approx}
u^h + \Pi_h \mathcal{T} \circ \mathcal{G}(u^h) = 0.
\ee
\el

\begin{Proof}
We organize the proof into two steps. Step 1 deals with the original problem \eqref{Eq_nonlinear}; and Step 2 deals with its {\tt lsfem} formulation.

\smallskip
{\bf Step 1.} Note that since $y_*$ solves \eqref{Eq_nonlinear}, it is a minimizer of the objective function $J$ given by \eqref{Def_J}. As a result, $y_*$ satisfies the following first-order necessary condition:
\be \label{Eq_1st_order_cond_v1}
\frac{\d}{\d \tau} J(y_*+\tau v; F, f, g) \big\vert_{\tau = 0} = 0, \qquad \forall \; v \in X.
\ee
By using \eqref{Def_J} in \eqref{Eq_1st_order_cond_v1}, we obtain the following integral equation to be satisfied by $y_*$:
\be \label{Eq_1st_order_cond}
\int_0^T \langle u' - F(u) - f, v' - \DF(u) v \rangle \, \d t + \langle u(0) - g, v(0) \rangle = 0,  \qquad \forall \; v \in X,
\ee
where $\DF$ denotes the Jacobian matrix of $F$.

In order to rewrite \eqref{Eq_1st_order_cond} into the form $u + \mathcal{T} \circ \mathcal{G}(u) = 0$,
let us introduce:
\be \label{Def_Q}
Q(u,v; f,g) =\! \int_0^T \langle u' -  F(u) - f,  \DF(u) v \rangle \, \d t  +\! \int_0^T \langle  F(u) + f, v' \rangle\, \d t + \langle g, v(0) \rangle,
\ee
which is defined for any $u,v$ in $X$, $f$ in $L^2(0, T; \mathbb{R}^d)$ and  $g$ in $\mathbb{R}^d$. Note that, with the above definition of $Q$ and the definition of the inner product on $X$ given by \eqref{Eq_X_inner},  \eqref{Eq_1st_order_cond} can be rewritten as
\be \label{Eq_1st_order_condition_rewritten}
\langle u,v \rangle_X - Q(u,v; f,g) = 0, \quad \forall \; v \in X.
\ee

To proceed further, we note that $Q(u,v; f,g)$ defined in \eqref{Def_Q} satisfies
\be \label{Eq_G_identity}
Q(u,v; f,g)  +  \int_0^T \langle \widetilde{f}(u), v' \rangle \,\d t + \langle \widetilde{g}(u), v(0) \rangle = 0, \qquad \forall \; u, v \in X,
\ee
where $(\widetilde{f}(u), \widetilde{g}(u))$ is an element in $L^2(0, T;  \mathbb{R}^d) \times \mathbb{R}^d$ given by \eqref{Eq_f_g_tilde}.

The above identity \eqref{Eq_G_identity} can be derived from integrating by parts the first term in the definition of $Q(u,v; f,g)$ given in \eqref{Def_Q}. Indeed, denoting $w = u' -  F(u) - f$ and by noting that
$\langle u' -  F(u) - f, \DF(u) v \rangle = \langle  w, \DF(u) v \rangle = ( \DF(u) v )\t w =
v\t (\DF(u))\t w$, we have
\bea \label{Eq_G_identity_justification}
\int_0^T \langle w, \DF(u) v \rangle\, \d t  &=  \int_0^T v\t (\DF(u))\t w \,\d t  \\
 & = - \int_0^T v\t \, \d\Big(\int_{t}^T (\DF(u))\t w \,\d s \Big) \\
& = \Big \langle
\int_{0}^T (\DF(u))\t w\,\d s, v(0) \Big \rangle  +  \int_0^T \Big \langle
\int_{t}^T (\DF(u))\t w \,\d s,  v' \Big \rangle \,\d t,
\eea
where we used integration by parts to obtain the last equality above. Now, replacing the first term on the RHS of  \eqref{Def_Q} with the RHS of  \eqref{Eq_G_identity_justification} and using the definition of $(\widetilde{f}(u), \widetilde{g}(u))$ given by \eqref{Eq_f_g_tilde}, we obtain \eqref{Eq_G_identity}.

Thanks to the identity \eqref{Eq_G_identity}, an element $u$ in $X$ satisfies \eqref{Eq_1st_order_condition_rewritten} if and only if
\be \label{Eq_1st_order_equiv_formulation}
\langle u,v \rangle_X = \int_0^T \langle -\widetilde{f}(u), v' \rangle\, \d t + \langle -\widetilde{g}(u), v(0) \rangle,  \quad \forall \; v \in X.
\ee
Recalling the equivalent characterization given in \eqref{Eq_T_equiv_formulation} of the linear operator $\mathcal{T}$ defined by \eqref{Eq_def_T}, we get from \eqref{Eq_1st_order_equiv_formulation} that
\bes
u = \mathcal{T}(-\widetilde{f}(u), -\widetilde{g}(u)), \text{ or equivalently, } u + \mathcal{T}(\widetilde{f}(u), \widetilde{g}(u)) = 0.
\ees
At the same time, by the definition of $\mathcal{G}$ in \eqref{Eq_def_G}, we have $(\widetilde{f}(u), \widetilde{g}(u)) = \mathcal{G}(u)$. Using this relation in $u + \mathcal{T}(\widetilde{f}(u), \widetilde{g}(u)) = 0$, we obtain the desired form $u + \mathcal{T} \circ \mathcal{G}(u) = 0$ for the first-order optimality condition \eqref{Eq_1st_order_cond}.

\smallskip
{\bf Step 2.} Now, we consider the {\tt lsfem} problem \eqref{Eq_opt_FEM} which aims to approximate the variational formulation \eqref{Eq_opt} of \eqref{Eq_nonlinear}. Note that any solution $y^h_*$ of \eqref{Eq_opt_FEM}, if exists, satisfies the analogue of \eqref{Eq_1st_order_cond} with $v$ therein restricted to $X^h$. That is
\be \label{Eq_1st_order_cond_FEM}
\int_0^T \langle (y^h_*)' - F(y^h_*) - f, (v^h)' - \DF(y^h_*) v^h \rangle\, \d t + \langle y^h_*(0) - g, v^h(0) \rangle = 0,  \; \forall \; v^h \in X^h.
\ee
We can then follow the same derivation of \eqref{Eq_1st_order_cond_operator_form} from \eqref{Eq_1st_order_cond} to obtain that $y^h_*$ is a solution of the following nonlinear problem defined on $X^h$:
\be \label{Eq_1st_order_cond_operator_form_Xh}
u^h + \mathcal{T}^h \circ \mathcal{G}(u^h) = 0, \quad u^h \in X^h,
\ee
where $\mathcal{G}$ is the same as defined in \eqref{Eq_def_G}, and  $\mathcal{T}^h \colon L^2(0, T; \mathbb{R}^d) \times \mathbb{R}^d \rightarrow X^h$ is defined by
\bea \label{Eq_def_Th}
 \mathcal{T}^h(f,g)  &= u^h \quad \text{if and only if}  \\
& \langle u^h,  v^h \rangle_X = \int_0^T \langle f ,  (v^h)' \rangle \, \d t + \langle g, v^h(0) \rangle,  \qquad \forall \; v^h \in X^h.
\eea

Comparing \eqref{Eq_1st_order_cond_operator_form_Xh} with \eqref{Eq_1st_order_cond_operator_form_approx}, it remains to show that
\be \label{Eq_proj_identity_T_appendix}
\mathcal{T}^h = \Pi_h \mathcal{T}.
\ee
To see this, for any $(f,g) \in L^2(0, T; \mathbb{R}^d) \times \mathbb{R}^d$, we get from \eqref{Eq_T_equiv_formulation} and \eqref{Eq_def_Th} that
\bes
\langle \mathcal{T}(f,g) - \mathcal{T}^h(f,g),v \rangle_X = 0, \quad \forall\, v \in X^h.
\ees
Namely, $\mathcal{T}(f,g) - \mathcal{T}^h(f,g)$ belongs to the orthogonal complement of $X^h$; that is $\mathcal{T}^h = \Pi_h \mathcal{T}$. The proof is complete.
\qed
\end{Proof}

Thanks to Lemma~\ref{Lem:operator_form}, we have thus reformulated the first-order optimality condition associated with each of the minimization problems \eqref{Eq_opt} and \eqref{Eq_opt_FEM} into the desired form given by \eqref{Eq_1st_order_cond_operator_form} and \eqref{Eq_1st_order_cond_operator_form_approx}, respectively.

Note that \eqref{Eq_1st_order_cond_operator_form} and \eqref{Eq_1st_order_cond_operator_form_approx} fit into the abstract formulation of \cite[Theorem 3.3, p.307]{girault2012finite} (see also \cite[Theorem 8.1]{Bochev2009}).
For a given nonsingular solution $y_*$ of  \eqref{Eq_1st_order_cond_operator_form},  \cite[Theorem 3.3, p.307]{girault2012finite} delineates conditions to ensure the existence of a solution for \eqref{Eq_1st_order_cond_operator_form_approx} for all sufficiently small $h$ that converges to $y_*$.\footnote{The formulation of \cite[Theorem 3.3, p.307]{girault2012finite} concerns actually a parameterized family of \eqref{Eq_1st_order_cond_operator_form} in which the nonlinearity $\mathcal{G}$ depends on an additional scalar parameter $\lambda$. Since there is no such a parameter in our setting, it can be viewed as a special case for which $\lambda$ is taken to be a constant here.}

A solution $y_*$ to \eqref{Eq_1st_order_cond_operator_form} is called nonsingular if the linear operator $\mathrm{Id}_X + \mathcal{T} \circ {\rm D} \mathcal{G}(y_*) \colon X \rightarrow X$ is invertible with bounded inverse; namely,
\be  \label{Eq_nonsingular_assumption}
(\mathrm{Id}_X + \mathcal{T} \circ {\rm D} \mathcal{G}(y_*))^{-1} \in L(X,X),
\ee
where $\mathrm{Id}_X$ denotes the identity map on $X$, and $L(X,X)$ the set of all bounded linear operators on $X$. As will be shown below that, for the problem at hand, the condition \eqref{Eq_nonsingular_assumption} is ensured by the smallness assumption on the operator norm of the Jacobian matrix $\DF(y_*(t))$ for all $t$ in $[0, T]$.

\begin{Proof}[\bf Proof of Theorem~\ref{Thm_convergence}]

From what precedes, to prove Theorem~\ref{Thm_convergence}, it remains to check that
\bi
\item[(i)] all the conditions required in \cite[Theorem 3.3, p.307]{girault2012finite} are satisfied, ensuring thus the existence of a solution $y^h_*$ of \eqref{Eq_1st_order_cond_operator_form_approx} for all sufficiently small $h$; and that

\item[(ii)] the solution $y^h_*$ of the first-order optimality condition \eqref{Eq_1st_order_cond_operator_form_approx} obtained from (i) is indeed a minimizer of \eqref{Eq_opt_FEM}.
\ei

We proceed in two steps below.

\medskip
\noindent{\bf Step 1.} Introducing the space $Y:= L^2(0, T; \mathbb{R}^d) \times \mathbb{R}^d$, the conditions required in \cite[Theorem 3.3, p.307]{girault2012finite} are:
\bi
\item[\bf{(C1)}] $y_*$ is a nonsingular solution to \eqref{Eq_1st_order_cond_operator_form} in the sense that \eqref{Eq_nonsingular_assumption} holds.
\item[{\bf(C2)}] $\mathcal{G} \colon X \rightarrow Y$ is $C^2$ smooth, and ${\rm D}^2 \mathcal{G}$ is bounded on all bounded subsets of $X$.
\item[{\bf(C3)}] There exists a subspace $Z$ of $Y$, with continuous embedding, such that ${\rm D} \mathcal{G}(u) \in L(X, Z), \quad \forall \; u \in X$.
\item[{\bf(C4)}] $\lim_{h\rightarrow 0} \|(\mathcal{T} - \mathcal{T}^h) (f,g)\|_{X} = 0,  \quad \forall \; (f,g) \in Y$.
\item[{\bf(C5)}] $\lim_{h\rightarrow 0} \|\mathcal{T} - \mathcal{T}^h\|_{L(Z,X)} = 0$.
\ei

\smallskip
\noindent{\bf Verification of {\bf(C1)}}. To guarantee \eqref{Eq_nonsingular_assumption}, it suffices to show that
\be \label{Eq_cond_C1_equiv}
\|\mathcal{T} \circ {\rm D} \mathcal{G}(y_*)\|_{L(X,X)} < 1.
\ee
To this end, for any given $v$ in $X$, let us denote $w = \mathcal{T} \circ {\rm D} \mathcal{G}(y_*) \, v$. By a direct calculation using the definition of $\mathcal{T}$ and $\mathcal{G}$ given respectively in \eqref{Eq_def_T} and \eqref{Eq_def_G}, we get
\beas
w (t) & = -\int_0^t \DF(y_*(s)) v(s)\,\d s \\
& \quad - \int_0^t \int_s^T [\DF(y_*(\tau))]^\top \Big(v'(\tau) - \DF(y_*(\tau)) v(\tau) \Big) \, \d \tau \,\d s \\
 & \quad - \int_0^T [\DF(y_*(s))]^\top \Big(v'(s) - \DF(y_*(s)) v(s) \Big)\, \d s, \quad t \in [0, T].
\eeas
Introducing
\be
M:= \sup_{s \in [0,T]} \| \DF(y_*(s))\|_{\mathrm{op}}^2,
\ee
we obtain by a direct estimation based on the H\"older's inequality that
\bea \label{Eq_est_TGV1}
\langle w(0), w(0) \rangle &\le 2 T M \big(\|v'\|^2_{L^2(0,T; \mathbb{R}^d)} + M \|v\|^2_{L^2(0,T; \mathbb{R}^d)} \big)  \\
& \le 2 T M \big(1 +  \widetilde{\mathfrak{C}}^2 M \big) \|v\|^2_X,
\eea
and that
\be \label{Eq_est_TGV2}
\int_0^T \langle w'(t), w'(t) \rangle\, \d t \le M (2\widetilde{\mathfrak{C}}^2 + 4T^2 + 4T^2 \widetilde{\mathfrak{C}}^2 M) \|v\|^2_X.
\ee
Recalling that $w = \mathcal{T} \circ {\rm D} \mathcal{G}(y_*) \, v$, we get from \eqref{Eq_est_TGV1} and \eqref{Eq_est_TGV2} that
\be
\|\mathcal{T} \circ {\rm D} \mathcal{G}(y_*) \, v\|_X^2 \le 2 M(2T^2 + T + \widetilde{\mathfrak{C}}^2 + T \widetilde{\mathfrak{C}}^2 (1 + 2T ) M)\|v\|^2_X.
\ee
Since $T$ is fixed, we get $\|\mathcal{T} \circ {\rm D} \mathcal{G}(y_*) \, v\|_X < \|v\|^2_X$ for all $v$ in $X$ when $M$  satisfies
\bes
2 M(2T^2 + T + \widetilde{\mathfrak{C}}^2 +  T \widetilde{\mathfrak{C}}^2 (1+2T) M) < 1.
\ees
That is
\bea \label{eq_M_cond}
M &= \sup_{s \in [0,T]} \| \DF(y_*(s))\|^2_{\mathrm{op}} \\
& < \frac{1}{2T^2 + T + \widetilde{\mathfrak{C}}^2 + \sqrt{(2T^2 + T + \widetilde{\mathfrak{C}}^2)^2 + 2 T \widetilde{\mathfrak{C}}^2(1 + 2T)}}.
\eea
We have thus verified \eqref{Eq_cond_C1_equiv} when $\sup_{s \in [0,T]} \| \DF(y_*(s))\|_{\mathrm{op}}$ is small such that \eqref{eq_M_cond} holds. Consequently, Condition {\bf(C1)} holds under this smallness assumption on $\DF$.

\smallskip
\noindent{\bf Verification of {\bf(C2)}}. This condition can be checked by a long but straightforward calculation using the explicit form of $\mathcal{G}$ given by \eqref{Eq_def_G}--\eqref{Eq_f_g_tilde} and the assumption that $F$ is $C^3$ smooth.

\smallskip
\noindent{\bf Verification of {\bf(C3)}}. We take $Z:= H^1(0, T; \mathbb{R}^d) \times \mathbb{R}^d$. Note that $Z$ is compactly embedded into $Y$. Condition {\bf(C3)} follows then from a direct but lengthy calculation based again on the explicit form of $\mathcal{G}$ given by \eqref{Eq_def_G}--\eqref{Eq_f_g_tilde}. It suffices to assume $F$ to be $C^2$ smooth to check this condition.

\smallskip
\noindent{\bf Verification of {\bf(C4)}}. Recall from \eqref{Eq_proj_identity_T_appendix} that $\mathcal{T}^h = \Pi_h \mathcal{T}$.
Condition {\bf(C4)} follows immediately because $\|(\mathrm{Id}_X - \Pi_h) u\|_X$ converges to zero as $h$ goes to zero for all $u\in X$ as a property of the finite element subspaces $X^h$.

\smallskip
\noindent{\bf Verification of {\bf(C5)}}. As pointed out in \cite[Theorem 3.3, p.307]{girault2012finite}, Condition {\bf(C5)} is a consequence of Condition {\bf(C4)} (and the uniform boundedness theorem) when $Z$ is compactly embedded into $Y$, which is the case here for $Z= H^1(0, T; \mathbb{R}^d) \times \mathbb{R}^d$. See also \cite[Lemma 8.7]{Bochev2009}.

All the conditions in \cite[Theorem 3.3, p.307]{girault2012finite} are thus verified. It follows then from this theorem that for any given neighborhood $\mathcal{O}$ of $y_*$, the problem \eqref{Eq_opt_FEM} has a unique solution $y_*^h$ in $\mathcal{O}$ for all sufficiently small $h$; and the convergence result \eqref{Eq_conv_Sect2} holds.

\medskip
\noindent{\bf Step 2.} It remains to show that $y_*^h$ obtained from Step 1 above is indeed a minimizer of \eqref{Eq_opt_FEM}. For this purpose, it suffices to show that
\be \label{Eq_2nd_optimality_cond}
\frac{\d^2}{\d \tau^2} J(y_*^h + \tau v; F, f, g) \vert_{\tau = 0} > 0, \quad \forall \; v \in X^h \backslash \{ 0 \}.
\ee
Note that
\beas
\frac{\d^2}{\d \tau^2} & J(y_*^h + \tau v; F, f, g) \big\vert_{\tau = 0} \\
& = \int_0^T \langle -\DtF(y_*^h(t)) (v(t),v(t)), (y_*^h)'(t) - F(y_*^h(t)) -f(t) \rangle \, \d t \\
& \quad + \int_0^T \langle v'(t) -\DF(y_*^h(t)) v(t), v'(t) - \DF(y_*^h(t)) v(t) \rangle\,  \d t + \langle v(0), v(0) \rangle,
\eeas
where $\DtF(y_*^h(t))$ denotes the Hessian of $F$ evaluated at $y_*^h(t)$, which is a bilinear function mapping $\mathbb{R}^d \times \mathbb{R}^d$ to $\mathbb{R}^d$.

A straightforward estimation leads then to
\bea \label{Eq_2nd_derivative_est_J}
\frac{\d^2}{\d \tau^2} & J(y_*^h + \tau v; F, f, g) \big\vert_{\tau = 0} \ge \|v\|^2_X  \Bigg( 1 - 2 \widetilde{\mathfrak{C}} \sup_{t \in [0,T]} \| \DF(y^h_*(t))\|_{\mathrm{op}} \\
& - \mathfrak{C} \widetilde{\mathfrak{C}} \Big(\sup_{t \in [0,T]} \vertiii{\DtF(y^h_*(t))}  \Big) \|(y_*^h)' - F(y_*^h) -f\|_{L^2(0,T; \mathbb{R}^d)}\Bigg),
\eea
where $\vertiii{\DtF(y_*(t))}$ denotes the operator norm of the bilinear map $\DtF(y_*(t))$, $\mathfrak{C}$ and $\widetilde{\mathfrak{C}}$ are the embedding constants defined at the end of Section~\ref{Sect_theory_prelim}.

Since it has been shown in Step 1 that $y_*^h$ converges in $X$-norm to $y_*$ (cf.~\eqref{Eq_conv_Sect2}) and $X$ is continuously embedded into $C([0,T]; \mathbb{R}^d)$, we get
\be \label{Eq_conv_in_Linfty}
\lim_{h \rightarrow 0} \max_{t\in [0,T]} \|y_*(t) - y_*^h(t)\| = 0.
\ee
It follows that
\be \label{Eq_conv_DF}
\lim_{h \rightarrow 0} \sup_{t \in [0,T]} \| \DF (y^h_*(t))\|_{\mathrm{op}} = \sup_{t \in [0,T]} \| \DF (y_*(t))\|_{\mathrm{op}}.
\ee
The convergence results \eqref{Eq_conv_Sect2} and \eqref{Eq_conv_in_Linfty} together with the smoothness of $F$ also imply that
\be \label{Eq_conv_VF}
\lim_{h \rightarrow 0} \|(y_*^h)' - F(y_*^h) -f\|_{L^2(0,T; \mathbb{R}^d)} =  \|(y_*)' - F(y_*) -f\|_{L^2(0,T; \mathbb{R}^d)} = 0,
\ee
where the second equality holds since $y_*$ is a strong solution of the IVP \eqref{Eq_nonlinear}.

Thanks also to \eqref{Eq_conv_in_Linfty}, we know that $\sup_{t \in [0,T]} \vertiii{{\rm D}^2F(y^h_*(t))}$ is uniformly bounded with respect to $h$. This uniform boundedness together with \eqref{Eq_conv_DF} and \eqref{Eq_conv_VF} implies that for all nonzero $v$ in $X^h$, the right-hand side of \eqref{Eq_2nd_derivative_est_J} is positive when $h$ is sufficiently small provided that
\be \label{Eq_DF_cond2}
\sup_{t \in [0,T]} \| \DF(y_*(t))\|_{\mathrm{op}} < \frac{1}{2 \widetilde{\mathfrak{C}}}.
\ee
We have thus verified \eqref{Eq_2nd_optimality_cond} under the condition \eqref{Eq_DF_cond2} by taking $h$ sufficiently small. The proof is now complete.
\qed
\end{Proof}

\section{Proofs of Lemma~\ref{Lem_errest_F=0}} \label{Sect_Linear_case_estimation}
 Note that \eqref{Eq_soln_linear_LSFEM} always has a solution since $X^h$ is finite dimensional and the objective function is bounded below by zero. The fact $\widetilde{y}^h_* = \Pi_h \widetilde{y}_*$ is the unique solution to \eqref{Eq_soln_linear_LSFEM} follows directly by inspecting the associated first-order optimality condition. This condition can be obtained from \eqref{Eq_1st_order_cond} by setting $F$ to zero and restricting $v$ to $X^h$, and it reads as follows
 \be \label{Eq_optimality_cond_linear_FEM}
\int_0^T \langle (\widetilde{y}_*^h)' - f, (v^h)'  \rangle \, \d t + \langle \widetilde{y}_*^h(0) - g, v^h(0) \rangle = 0,  \qquad \forall \; v^h \in X^h.
\ee
Using the expression of the solution $\widetilde{y}_*$ given by \eqref{Eq_soln_linear} and the definition of the inner product $\langle \cdot, \cdot \rangle_X$ given by \eqref{Eq_X_inner}, we can rewrite the above condition as
\be
\langle \widetilde{y}_*^h - \widetilde{y}_*,  v^h \rangle_X = 0, \qquad  \forall \; v^h \in X^h.
\ee
Hence, $\widetilde{y}^h_* - \widetilde{y}_*$ lives in the orthogonal complement of $X^h$. We get thus, $\widetilde{y}^h_* = \Pi_h \widetilde{y}_*$.

For the error estimate \eqref{Eq_error_estimate_trivial_case}, see e.g.,~\cite[Section 2.7.3]{Jiang1998} for a proof that relies on the classical Aubin-Nitsche trick. The proof presented therein deals with the special case $g=0$ and for state space dimension $d = 1$. For $d>1$, since the vector field is independent of the unknown variable, we can carry out the estimate component by component, which reduces the problem to the case $d=1$. The general case of $g \neq 0$ can be handled by considering $z = \widetilde{y} - g$.
\qed

\bibliographystyle{plain}
\bibliography{lsfem_arXiv}

\begin{thebibliography}{10}

\bibitem{ascher1994numerical}
U~M Ascher, R~M~M Mattheij, and R~D Russell.
\newblock {\em {Numerical Solution of Boundary Value Problems for Ordinary
  Differential Equations}}, volume~13 of {\em Classics in Applied Mathematics}.
\newblock SIAM, Philadelphia, PA, 1995.

\bibitem{Ascher1998}
U~M Ascher and L~R Petzold.
\newblock {\em Computer Methods for Ordinary Differential Equations and
  Differential-Algebraic Equations}.
\newblock SIAM, Philadelphia, PA, 1998.

\bibitem{aubin1967behavior}
J~P Aubin.
\newblock {Behavior of the error of the approximate solutions of boundary value
  problems for linear elliptic operators by {G}alerkin's and finite difference
  methods}.
\newblock {\em Annali della Scuola Normale Superiore di Pisa, Classe di
  Scienze}, 21:599--637, 1967.

\bibitem{Becker_al1981}
E~B Becker, G~F Carey, and J~T Oden.
\newblock {\em {Finite Elements, An Introduction: Volume I}}.
\newblock Prentice-Hall, Inc., Englewood Cliffs, NJ, 1981.

\bibitem{bellen2013numerical}
A~Bellen and M~Zennaro.
\newblock {\em {Numerical Methods for Delay Differential Equations}}.
\newblock Oxford University Press, Oxford, 2013.

\bibitem{Bochev2009}
P~B Bochev and M~D Gunzburger.
\newblock {\em Least-Squares Finite Element Methods}.
\newblock Springer, New York, 2009.

\bibitem{brezis_book}
H~Br{\'{e}}zis.
\newblock {\em Functional {Analysis, Sobolev Spaces and Partial Differential
  Equations}}.
\newblock Springer, New York, 2010.

\bibitem{ciarlet2002finite}
P~G Ciarlet.
\newblock {\em The finite element method for elliptic problems}, volume~40 of
  {\em Classics in Applied Mathematics}.
\newblock SIAM, Philadelphia, 2002.

\bibitem{Cockburn00}
B~Cockburn, G~E Karniadakis, and C-W Shu.
\newblock The development of discontinuous {G}alerkin methods.
\newblock In {\em {Discontinuous Galerkin Methods}}, volume~11 of {\em Lect.
  Notes Comput. Sci. Eng.}, pages 3--50. Springer, Berlin, 2000.

\bibitem{Davis1984}
M~E Davis.
\newblock {\em {Numerical Methods and Modeling for Chemical Engineers}}.
\newblock John Wiley \& Sons, Inc., New York, 1984.

\bibitem{deBoor2001}
C~de~Boor.
\newblock {\em A Practical Guide to Splines}.
\newblock Springer-Verlag, New York, NY, 1st edition, 2001.

\bibitem{eilertsen2020quasi}
J~Eilertsen and S~Schnell.
\newblock The quasi-steady-state approximations revisited: Timescales, small
  parameters, singularities, and normal forms in enzyme kinetics.
\newblock {\em Mathematical Biosciences}, 325:108339, 2020.

\bibitem{Evans10}
L~C Evans.
\newblock {\em Partial {D}ifferential {E}quations}, volume~19 of {\em Graduate
  Studies in Mathematics}.
\newblock American Mathematical Society, Providence, RI, 2010.

\bibitem{Fairweather89}
G~Fairweather and D~Meade.
\newblock A survey of spline collocation methods for the numerical solution of
  differential equations.
\newblock In {\em {Mathematics for Large Scale Computing}}, volume 120 of {\em
  Lecture Notes in Pure and Appl. Math.}, pages 297--341. Dekker, New York,
  1989.

\bibitem{girault2012finite}
V~Girault and P-A Raviart.
\newblock {\em {Finite Element Methods for Navier-Stokes Equations: Theory and
  Algorithms}}.
\newblock Springer, New York, 1986.

\bibitem{Golub1967}
G~H Golub and J~H Welsch.
\newblock Calculation of {G}auss quadrature rules.
\newblock Technical report, Stanford, CA, USA, 1967.

\bibitem{griewank2008evaluating}
A~Griewank and A~Walther.
\newblock {\em Evaluating derivatives: principles and techniques of algorithmic
  differentiation}, volume 105.
\newblock SIAM, 2008.

\bibitem{gunther2021spline}
S~G{\"u}nther, W~Pazner, and D~Qi.
\newblock Spline parameterization of neural network controls for deep learning.
\newblock {\em arXiv preprint arXiv:2103.00301}, 2021.

\bibitem{haber2017stable}
E~Haber and L~Ruthotto.
\newblock Stable architectures for deep neural networks.
\newblock {\em Inverse problems}, 34(1):014004, 2017.

\bibitem{hairer2006geometric}
E~Hairer, C~Lubich, and G~Wanner.
\newblock {\em {Geometric Numerical Integration: Structure-Preserving
  Algorithms for Ordinary Differential Equations}}, volume~31 of {\em Springer
  Series in Computational Mathematics}.
\newblock Springer, New York, 2nd edition, 2006.

\bibitem{Hairer1993}
G~Hairer, S~P N{\o}rsett, and E~Wanner.
\newblock {\em Solving Ordinary Differential Equations I: Nonstiff Problems}.
\newblock Springer, New York, 2nd edition, 1993.

\bibitem{Hairer1996}
G~Hairer and E~Wanner.
\newblock {\em Solving Ordinary Differential Equations II: Stiff and
  Differential-Algebraic Problems}.
\newblock Springer, New York, 2nd edition, 1996.

\bibitem{horst2013handbook}
R~Horst and P~M Pardalos.
\newblock {\em Handbook of Global Optimization}, volume~2.
\newblock Springer, 2013.

\bibitem{reis2014surveys}
A~Ilchmann and T~Reis.
\newblock {\em Surveys in Differential-algebraic Equations II}.
\newblock Springer, 2014.

\bibitem{Jiang1998}
B~Jiang.
\newblock {\em The Least-Squares Finite Element Method: Theory and Applications
  in Computational Fluid Dynamics and Electromagnetics}.
\newblock Springer, New York, 1998.

\bibitem{kreyszig11}
E~Kreyszig, H~Kreyszig, and E~J Norminton.
\newblock {\em Advanced Engineering Mathematics}.
\newblock Wiley, Hoboken, NJ, tenth edition, 2011.

\bibitem{lambert1991numerical}
J~D Lambert.
\newblock {\em {Numerical Methods for Ordinary Differential Systems: The
  Initial Value Problem}}.
\newblock Wiley New York, 1991.

\bibitem{nitsche1968kriterium}
J~Nitsche.
\newblock {Ein Kriterium f{\"u}r die Quasi-Optimalit{\"a}t des Ritzschen
  Verfahrens}.
\newblock {\em Numerische Mathematik}, 11:346--348, 1968.

\bibitem{Rentrop1980}
P~Rentrop.
\newblock An algorithm for the computation of the exponential spline.
\newblock {\em Numerische Mathematik}, 35(1):81--93, 1980.

\bibitem{Royden2010real}
H~L Royden and P~M Fitzpatrick.
\newblock {\em {Real Analysis}}.
\newblock Pearson, Boston, MA, 4th edition, 2010.

\bibitem{Schumaker2015}
L~L Schumaker.
\newblock {\em Spline Functions: Computational Methods}.
\newblock SIAM, Philadelphia, PA, 2015.

\bibitem{segel1988validity}
L~A Segel.
\newblock On the validity of the steady state assumption of enzyme kinetics.
\newblock {\em Bulletin of Mathematical Biology}, 50(6):579--593, 1988.

\end{thebibliography}

\end{document}